\documentclass[11pt]{amsart}
\usepackage[a4paper,margin=3cm]{geometry}
\usepackage[T1]{fontenc}
\usepackage{amsmath,amssymb,mathtools,microtype}
\usepackage{algorithm}
\usepackage{algpseudocode}
\usepackage[most]{tcolorbox}
\usepackage{pgfplots}
\usepgfplotslibrary{groupplots}
\pgfplotsset{compat=1.18}
\usepackage[hidelinks]{hyperref}
\definecolor{algogray}{gray}{0.96}
\tcbset{algobox/.style={
  colback=algogray,
  colframe=gray!50,
  boxrule=0.7pt,
  arc=2pt,
  left=7pt,
  right=7pt,
  top=7pt,
  bottom=7pt,
  enhanced
}}
\newtheorem{theorem}{Theorem}[section]
\newtheorem{proposition}[theorem]{Proposition}
\newtheorem{corollary}[theorem]{Corollary}
\theoremstyle{remark}
\newtheorem{remark}[theorem]{Remark}
\newtheorem{example}[theorem]{Example}
\newcommand{\Pp}{\mathcal P}
\newcommand{\op}[1]{\mathsf{#1}}
\newcommand{\Id}{\op{I}}
\newcommand{\eqtag}[1]{%
  \leqno\hbox{\normalfont\refstepcounter{equation}\label{#1}(\theequation)}%
}
\AtBeginDocument{%
  \setlength{\abovedisplayskip}{10pt plus 2pt minus 2pt}%
  \setlength{\belowdisplayskip}{10pt plus 2pt minus 2pt}%
  \setlength{\abovedisplayshortskip}{5pt plus 2pt minus 1pt}%
  \setlength{\belowdisplayshortskip}{7pt plus 2pt minus 1pt}%
}

\title[Another application of the alternating map]
{Another application of the alternating map: time-and-band limiting}
\author{L. Tertuliano da Silva}
\address{Department of Mathematics and Computer Science,
S\~ao Paulo State University (UNESP), Presidente Prudente, SP, Brazil}
\email{lucas.tertuliano@unesp.br}
\date{}
\hypersetup{pdftitle={Another application of the alternating map: time-and-band limiting},
pdfauthor={Lucas Tertuliano da Silva}}

\begin{document}

\begin{abstract}
Castillo's dual equations for the alternating map reduce polynomial
time-and-band limiting to two scalar problems. We determine the
concentration multiplicities at a single cutoff, with sharp finite
support and moment requirements. When the scalar problems admit
irreducible tridiagonal commutants, we classify every real symmetric
pentadiagonal commutant of the simple concentration matrix and
construct one with simple spectrum. Its half-bandwidth is minimal
for $N\ge3$ in polynomial degree order.
Laguerre, Jacobi and Meixner input give applications; the constructions
of Castro and Gr\"unbaum follow by specialization.
We give spectral-gap and rounding bounds, and certify small
concentration values and deficits for noninteger Jacobi parameters.
\end{abstract}

\maketitle

\section{Classicality for the alternating map}

We use Castillo's classification and dual reconstruction for the alternating map
\cite[Theorem~5.1(ii) and Section~8]{Admissible}.
The dual pullback transfers orthogonality, concentration forms and
operator symmetries from the input functional.
For related polynomial mappings developed by Petronilho and his students, see
\cite[Theorems~2.1 and~4.2]{PolynomialMappings}.
The measure formula there assumes compact input support.

Let $h\in\mathbb C\setminus\{0\}$ and
$U\subseteq\mathbb C$ satisfy
$$
U\pm h/2\subseteq U,\quad X(s)=e^{\pi i s/h}\ (s\in U),\quad \#X(U)=\infty.
$$
For example, $U=h\mathbb Q$ has infinite image.
Let $\Pp$ be the space of polynomials with complex coefficients and
$\Pp_n=\{p\in\Pp:\deg p\le n\}$, with $\Pp_{-1}=\{0\}$.
Write $\Pp'$ for its linear dual.
Since $X(s\pm h/2)=\pm iX(s)$, the divided-difference and averaging operators are
$$
(\op{D}p)(x)=\frac{p(ix)-p(-ix)}{2ix},\quad
(\op{S}p)(x)=\frac{p(ix)+p(-ix)}2.
$$
The first quotient extends polynomially through $x=0$.
In Maroni's dual notation \cite[Section~1.1]{Maroni1993}, the transpose of
an operator $\op{A}$ and multiplication of a functional are defined by
$$
\langle\op{A}'\mathbf u,p\rangle=\langle\mathbf u,\op{A}p\rangle,
\quad
\langle q\mathbf u,p\rangle=\langle\mathbf u,qp\rangle.
$$
For $\mathbf a\in\Pp'$, define $\op{J}_{\mathbf a}:\Pp\to\Pp'$ by
$$
\op{J}_{\mathbf a}f=f\mathbf a,\quad
\langle\op{J}_{\mathbf a}f,g\rangle=\langle\mathbf a,fg\rangle.
$$
An operator $\op{A}$ is symmetric for $\mathbf a$ when
$\op{J}_{\mathbf a}\op{A}=\op{A}'\op{J}_{\mathbf a}$.
All transposes refer to this bilinear pairing.
A functional is regular when the form is nondegenerate on each $\Pp_n$;
equivalently, it has a monic orthogonal sequence with nonzero squared
norms in every degree. A regular functional is classical for this
alternating map when there exist polynomials $\phi,\psi$, not both zero,
satisfying \cite[Definition~1.2]{Manufacture}
$$
-\op{D}'(\phi\mathbf u)=\op{S}'(\psi\mathbf u),\quad
\deg\phi\le2,\quad\deg\psi\le1.
$$
For real $\tau$, we use $\phi=0$, $\psi=x-\tau$. The equation becomes
$\op{S}'((x-\tau)\mathbf u)=0$.

\begin{samepage}
Let $\nu$ be a finite positive Borel measure on $(0,b)$,
$0<b\le\infty$, with infinite support and all moments finite, and let
$(R_n)_{n\ge0}$ be the monic real orthogonal polynomials of
$\langle\mathbf v,p\rangle=\int_{(0,b)}p(t)\,d\nu(t)$.

\end{samepage}
Set, for $n\ge0$,
$$
\begin{aligned}
S_n(t)&=\frac{R_{n+1}(t)-\rho_nR_n(t)}t,\quad
\rho_n=\frac{R_{n+1}(0)}{R_n(0)},\\[7pt]
h_n&=\langle\mathbf v,R_n^2\rangle,\\[7pt]
\widehat R_n&=\frac{R_n}{\sqrt{h_n}},\quad
\widehat S_n=\frac{S_n}{\sqrt{-\rho_nh_n}}.
\end{aligned}
\eqtag{eq:input}
$$
Here $R_n(0)\ne0$: otherwise $R_n(t)=tQ(t)$ and
$0=\langle\mathbf v,R_nQ\rangle
=\int_{(0,b)}tQ(t)^2\,d\nu(t)>0$.
The sequence $(S_n)$ is orthogonal for $x\mathbf v$, with squared norms
$-\rho_nh_n>0$.

For $\tau\in\mathbb R$, put
$\langle\mathbf v_\tau,p\rangle=\langle\mathbf v,p(x+\tau^2)\rangle$.
Every polynomial $F(x)=E(x^2)+xO(x^2)$ has the unique form
$F(x)=A(x^2)+(x-\tau)B(x^2)$, with $B=O$ and $A=E+\tau O$.
Following
\cite[Theorem~8.6]{Admissible} and \cite[Proposition~1.3]{Manufacture}
we define $\mathbf u_\tau$ by
$$
\langle\mathbf u_\tau,p(x^2)\rangle=\langle\mathbf v_\tau,p\rangle,
\quad
\op{S}'((x-\tau)\mathbf u_\tau)=0.
\eqtag{eq:functionals}
$$
\begin{samepage}
Write $y=x^2-\tau^2$ and
$E_\tau=\{x\in\mathbb R:0<x^2-\tau^2<b\}$.
Since $\op{S}(\Pp)$ is the space of even polynomials, the alternating
condition in \eqref{eq:functionals} means
$\langle\mathbf u_\tau,(x-\tau)p(x^2)\rangle=0$ for every $p\in\Pp$.
In particular,
$\langle\mathbf u_\tau,x^{2k+1}\rangle
=\tau\langle\mathbf u_\tau,x^{2k}\rangle$ for every $k\ge0$,
so the even moments determine the odd moments.
Let $\op{E}_\tau F=f$ and $\op{O}_\tau F=g$ in the unique decomposition
$F(x)=f(y)+(x-\tau)g(y)$. Then
$$
\mathbf u_\tau=\op{E}_\tau'\mathbf v.
$$
This proves uniqueness directly in $\Pp'$.

\end{samepage}

\begin{theorem}\label{thm:deformation}
The functional $\mathbf u_\tau$ is positive definite and classical for
the alternating map, with structural pair $\phi=0$, $\psi=x-\tau$.
With $r=\sqrt{t+\tau^2}$,
a representing measure is determined, for bounded Borel $F$, by
$$
\int_{E_\tau}F(x)\,d\mu_\tau(x)
=\int_{(0,b)}
\left(\frac{r+\tau}{2r}F(r)+\frac{r-\tau}{2r}F(-r)\right)d\nu(t).
\eqtag{eq:measure}
$$
Its orthonormal polynomials are
$$
\begin{aligned}
\widehat P_{2n}^{(\tau)}(x)&=\widehat R_n(y),\\[7pt]
\widehat P_{2n+1}^{(\tau)}(x)&=(x-\tau)\widehat S_n(y).
\end{aligned}
\eqtag{eq:deformed-polynomials}
$$
The operator $\op{U}_\tau$ defined by
$$
\bigl(\op{U}_\tau(f,g)\bigr)(x)=f(y)+(x-\tau)g(y)
\eqtag{eq:unitary}
$$
is unitary from
$L^2((0,b),d\nu(t))\oplus L^2((0,b),t\,d\nu(t))$
onto $L^2(E_\tau,d\mu_\tau(x))$
and intertwines multiplication by a bounded measurable function $\chi$ on
$(0,b)$ with multiplication by $\chi(y)$.

For an integer $N\ge1$ and measurable $\chi:(0,b)\to\mathbb R$
with $0\le\chi\le1$, define
$$
\begin{gathered}
\langle\mathbf v_\chi,p\rangle
=\int_{(0,b)}p(t)\chi(t)\,d\nu(t),\\[7pt]
\langle\mathbf u_{\tau,\chi},p\rangle
=\int_{E_\tau}p(x)\chi(x^2-\tau^2)\,d\mu_\tau(x),\\[7pt]
(\op{M}_\tau(\chi))_{ij}
=\langle\mathbf u_{\tau,\chi},
\widehat P_i^{(\tau)}\widehat P_j^{(\tau)}\rangle,\quad 0\le i,j<N.
\end{gathered}
$$
Then
$$
\op{M}_\tau(\chi)=\op{M}_0(\chi),\quad
\op{Q} \op{M}_\tau(\chi)\op{Q}^{\mathsf T}
=\op{G}_R(\chi)\oplus \op{G}_S(\chi),
\eqtag{eq:isospectral}
$$
where $\op{Q}$ orders the even indices increasingly, followed by
the odd indices increasingly. The block sizes are
$N_R=\lceil N/2\rceil$ and $N_S=\lfloor N/2\rfloor$, with entries
$$
\begin{aligned}
(\op{G}_R)_{mn}&=\langle\mathbf v_\chi,\widehat R_m\widehat R_n\rangle,
\quad 0\le m,n<N_R,\\[7pt]
(\op{G}_S)_{mn}&=\langle x\mathbf v_\chi,\widehat S_m\widehat S_n\rangle,
\quad 0\le m,n<N_S.
\end{aligned}
$$
\end{theorem}

We suppress the argument $\chi$ in $\op{G}_R,\op{G}_S$ when the mask is fixed.

\begin{proof}
For $f_\ell,g_\ell,q\in\Pp$, put
$F_\ell(x)=f_\ell(y)+(x-\tau)g_\ell(y)$, $\ell=1,2$.
The alternating condition annihilates $(x-\tau)h(y)$ for every
$h\in\Pp$. Using $(x-\tau)^2=y-2\tau(x-\tau)$ gives
$$
\begin{aligned}
\langle\mathbf u_\tau,q(y)F_1F_2\rangle
&=\langle\mathbf v,qf_1f_2+xqg_1g_2\rangle,\\[7pt]
\op{U}_\tau'\op{J}_{q(y)\mathbf u_\tau}\op{U}_\tau
&=\op{J}_{q\mathbf v}\oplus\op{J}_{xq\mathbf v}.
\end{aligned}
\eqtag{eq:dual-pullback}
$$
Here $\op{U}_\tau'$ acts from $\Pp'$ to
$(\Pp\oplus\Pp)'=\Pp'\oplus\Pp'$.
Taking $q=1$ proves the polynomial orthogonality in
\eqref{eq:deformed-polynomials}.
The branch weights in \eqref{eq:measure} are positive because
$r>|\tau|$, sum to one
and have weighted mean $\tau$ and variance $t$. Hence \eqref{eq:measure} satisfies
\eqref{eq:functionals}. All moments are finite since
$\int_{E_\tau}|x|^k\,d\mu_\tau(x)
=\int_{(0,b)}(t+\tau^2)^{k/2}\,d\nu(t)$.
The same branch calculation gives, for measurable complex-valued $f,g$,
$$
\int_{E_\tau}|f(y)+(x-\tau)g(y)|^2\,d\mu_\tau(x)
=\int_{(0,b)}\bigl(|f(t)|^2+t|g(t)|^2\bigr)\,d\nu(t)
\eqtag{eq:norm}
$$
with both integrals allowed to be infinite.
For a nonzero polynomial $F$, at least one of its components $f,g$ is
nonzero. Infinite support makes its squared norm strictly positive.
Thus $\mathbf u_\tau$ is regular and \eqref{eq:functionals} gives its
alternating classicality. Equation~\eqref{eq:norm} proves the isometry.
For $F\in L^2(E_\tau,d\mu_\tau(x))$ its inverse is given
$\nu$-almost everywhere by
$$
g(t)=\frac{F(r)-F(-r)}{2r},\quad
f(t)=\frac{F(r)+F(-r)}2+\tau g(t).
$$
The positive branch weights make $F(\pm r)$ well defined
$\nu$-almost everywhere. Substitution gives $\op{U}_\tau(f,g)=F$,
and \eqref{eq:norm} puts both components in the required spaces.
For the mask $\chi$, the branch measure gives the dual identities
$$
\begin{aligned}
\mathbf u_{\tau,\chi}&=\op{E}_\tau'\mathbf v_\chi,\\[7pt]
\op{U}_\tau'\op{J}_{\mathbf u_{\tau,\chi}}\op{U}_\tau
&=\op{J}_{\mathbf v_\chi}\oplus\op{J}_{x\mathbf v_\chi}.
\end{aligned}
\eqtag{eq:dual-mask}
$$
These identities hold for the constructed measures, without a moment
determinacy assumption. Multiplication by $\chi(y)$ acts componentwise.
Taking the matrices of the two sides of \eqref{eq:dual-mask} in the
reconstructed orthonormal bases proves \eqref{eq:isospectral}.

\end{proof}

Put $r_R=N_R-1$, $r_S=N_S-1$, with empty
matrix blocks omitted. The map \eqref{eq:unitary} carries
$\Pp_{r_R}\oplus\Pp_{r_S}$ onto $\Pp_{N-1}$.
Let $\Pi_R,\Pi_S,\Pi_\tau$ be the orthogonal projections onto these
three spaces in their respective $L^2$ spaces, and let $\op{Mult}_a$
denote multiplication by a bounded measurable $a$.
Unitarity and the multiplier identity give
$$
\begin{aligned}
\op{U}_\tau^{-1}\Pi_\tau\op{U}_\tau
&=\Pi_R\oplus\Pi_S,\\[7pt]
\op{U}_\tau^{-1}\Pi_\tau\op{Mult}_{\chi(y)}\Pi_\tau\op{U}_\tau
&=(\Pi_R\op{Mult}_\chi\Pi_R)\oplus(\Pi_S\op{Mult}_\chi\Pi_S).
\end{aligned}
$$
For $N=1$, $\Pi_S=0$. These identities use the measure
$\mu_\tau$ constructed in \eqref{eq:measure}; no determinacy or
polynomial density assumption is needed.

By \eqref{eq:norm}, the concentration problem is
$$
\begin{aligned}
\lambda_{\max}(\op{M}_\tau(\chi))
&=\max_{0\ne F\in\Pp_{N-1}}
\frac{\int_{E_\tau}|F(x)|^2\chi(x^2-\tau^2)\,d\mu_\tau(x)}
     {\int_{E_\tau}|F(x)|^2\,d\mu_\tau(x)}\\[7pt]
&=\max_{0\ne(f,g)\in\Pp_{r_R}\oplus\Pp_{r_S}}
\frac{\int_{(0,b)}\chi(t)(|f(t)|^2+t|g(t)|^2)\,d\nu(t)}
     {\int_{(0,b)}(|f(t)|^2+t|g(t)|^2)\,d\nu(t)}.
\end{aligned}
\eqtag{eq:variational}
$$
The optimal value is
$\max\{\lambda_{\max}(\op{G}_R),\lambda_{\max}(\op{G}_S)\}$,
with empty blocks omitted.
Every maximizing pair reconstructs an extremal polynomial for each $\tau$, of degree
$\max\{2\deg f,2\deg g+1\}$, with $\deg0=-\infty$.

\begin{theorem}\label{thm:concentration-multiplicities}
Let $0<s<b$ and $\chi=\mathbf1_{(0,s)}$. Let $m_-$ and $m_+$ be the
numbers of distinct support points of the restrictions of $\nu$ to
$(0,s)$ and $[s,b)$, respectively, allowing the value $\infty$.
Every eigenvalue of $\op{M}_\tau(\chi)$ strictly between zero and one
is simple. The multiplicities at the endpoints are
$$
\begin{aligned}
\operatorname{mult}(0)
&=N-\min(N_R,m_-)-\min(N_S,m_-),\\[7pt]
\operatorname{mult}(1)
&=N-\min(N_R,m_+)-\min(N_S,m_+),
\end{aligned}
\eqtag{eq:concentration-multiplicities}
$$
where $\min(n,\infty)=n$. Thus $\op{M}_\tau(\chi)$ has simple spectrum
precisely when both $m_-,m_+\ge\lfloor N/2\rfloor$.
In particular, its spectrum is simple for every $N$ if both
restricted measures have infinite support.
\end{theorem}

\begin{proof}
Fix $0<\lambda<1$ and define
$$
\mathbf a_\lambda=\mathbf v_\chi-\lambda\mathbf v,\quad
\mathbf w_\lambda=(s-x)\mathbf a_\lambda.
\eqtag{eq:positive-dual-cutoff}
$$
The measure representing $\mathbf w_\lambda$ is
$(s-t)(\chi(t)-\lambda)\,d\nu(t)$. It is positive away from $s$;
removing a possible atom at $s$ leaves infinite support.
Thus $\mathbf w_\lambda$ is positive definite.
Let $F_k$ be its monic orthogonal polynomials and put
$c_k=\langle\mathbf a_\lambda,F_k\rangle$.

An eigenvector of $\op{G}_R$ with eigenvalue $\lambda$ represents a
nonzero real polynomial $f\in\Pp_{N_R-1}$ satisfying
$\langle\mathbf a_\lambda,fq\rangle=0$ for $q\in\Pp_{N_R-1}$.
Testing with $q=(s-x)h$ gives
$\langle\mathbf w_\lambda,fh\rangle=0$ for $h\in\Pp_{N_R-2}$.
Positivity forces $\deg f=N_R-1$, so $f$ is a multiple of $F_{N_R-1}$.
The same argument for $x\mathbf a_\lambda$ treats $\op{G}_S$.
Each block eigenspace therefore has dimension one.

Suppose $N\ge2$ and $\lambda$ belongs to both blocks. Write
$n=N_S$, so $N_R=n$ or $n+1$, and let $G_{n-1}$ be monic orthogonal
for $x\mathbf w_\lambda$. Applying \eqref{eq:input} to
$\mathbf w_\lambda$ gives
$$
xG_{n-1}=F_n-\rho F_{n-1},\quad
\rho=\frac{F_n(0)}{F_{n-1}(0)}\ne0.
$$
The two eigenconditions, tested with $q=1$, imply
$c_{N_R-1}=0$ and $c_n-\rho c_{n-1}=0$. Hence $c_{n-1}=c_n=0$.
Set $F_{-1}=0$, $C_0=0$, and write the recurrence for $F_k$ as
$$
xF_k=F_{k+1}+B_kF_k+C_kF_{k-1},\quad C_k>0\ (k\ge1).
$$
Pairing with $\mathbf a_\lambda$ gives
$$
\begin{aligned}
c_{k+1}&=(s-B_k)c_k-C_kc_{k-1}\quad(k\ge1),\\[7pt]
c_1&=(s-B_0)c_0-\langle\mathbf w_\lambda,1\rangle.
\end{aligned}
$$
Consecutive zeros propagate backwards to $c_0=c_1=0$, contradicting
$\langle\mathbf w_\lambda,1\rangle>0$.
The blocks therefore have no common eigenvalue in $(0,1)$.
For $N=1$ there is only a scalar block.

For a restriction supported at $m$ distinct points, its polynomial
Gram matrix of order $n$ has rank $\min(n,m)$: evaluation at the
support points has this rank. Infinite support gives rank $n$.
Multiplication of the measure by $t>0$ does not change its support.
Apply these facts to the blocks of $\op{M}_\tau$ and $\Id-\op{M}_\tau$
to obtain \eqref{eq:concentration-multiplicities}.
Since all remaining eigenvalues are simple, the displayed criterion follows.
\end{proof}

\begin{proposition}\label{prop:finite-input}
For a fixed $N\ge1$, the finite reconstruction and
Theorem~\ref{thm:concentration-multiplicities} require only a nonzero
finite positive measure $\nu$ on $(0,b)$ with at least
$N_R=\lceil N/2\rceil$ distinct support points and
$$
\int_{(0,b)}t^{N-1}\,d\nu(t)<\infty.
$$
Define the real orthonormal polynomials
$\widehat R_0,\ldots,\widehat R_{N_R-1}$ for $\nu$ and
$\widehat S_0,\ldots,\widehat S_{N_S-1}$ for $t\,d\nu$,
with positive leading coefficients.
Then \eqref{eq:measure}, \eqref{eq:deformed-polynomials},
\eqref{eq:norm} and \eqref{eq:isospectral} hold for these finite
systems. Both assumptions are necessary for a nondegenerate
concentration problem on all of $\Pp_{N-1}$.
\end{proposition}

\begin{proof}
The required scalar moments have maximal order
$\max(2N_R-2,2N_S-1)=N-1$; for $N=1$ omit the second term.
The support condition makes both finite polynomial inner products
nondegenerate. The branch calculation proves the finite reconstruction.
Only the multiplicity argument needs a finite-moment proof.
We choose all eigenpolynomials below real.

For a positive measure with at least $k$ support points, an interior
eigenpolynomial of its size-$k$ concentration compression has degree
$k-1$. Otherwise, testing with $(s-x)f$ gives
$\int(s-t)(\chi(t)-\lambda)f(t)^2\,d\nu(t)=0$.
It forces $f$ to vanish on at least $k-1$ points away from $s$,
although $\deg f\le k-2$.
This also proves simplicity within each scalar block.
Two consecutive compressions have no common interior eigenvalue:
an eigenvector of the larger one has nonzero last coordinate;
the residual of an eigenvector of the smaller one, extended by zero,
is supported only in that coordinate. Pairing the two vectors forces
the residual to vanish, contradicting the degree assertion.

Put $\mathbf a=\mathbf v_\chi-\lambda\mathbf v$, understood only on
the degrees for which its moments exist.
If $N_R=N_S=n$ and $f,g$ are common eigenpolynomials, choose $c$
so $h=f-cg$ has degree at most $n-2$.
For $q\in\Pp_{n-2}$, both eigenconditions give
$\langle x\mathbf a,hq\rangle=0$.
The consecutive $t\,d\nu$ compressions force $h=0$.
For $n=1$ this is immediate.
Consequently
$\langle\mathbf a,f^2\rangle=\langle x\mathbf a,f^2\rangle=0$.
Their difference with factor $s$ forces $f$ to vanish off $s$;
the first equality forces it to vanish at any cutoff atom as well.
This contradicts its positive norm.

If $N_R=n+1$, $N_S=n$, choose $c$ so $h=xg-cf\in\Pp_{n-1}$.
Now $\langle\mathbf a,hq\rangle=0$ for $q\in\Pp_{n-1}$.
The consecutive $\nu$ compressions force $h=0$.
After rescaling, $f=xg$; testing its eigencondition with $f$ and $g$
gives $\langle x^2\mathbf a,g^2\rangle
=\langle x\mathbf a,g^2\rangle=0$.
Since $t(t-s)(\chi(t)-\lambda)<0$ off $s$, these equalities again
contradict the positive norm, including at a cutoff atom.
The case $N=1$ is scalar. Every product used has degree at most $N-1$;
the endpoint ranks are unchanged.

Finally, $M$ input support points give $2M$ reconstructed points, and
$$
\int_{E_\tau}x^{2N-2}\,d\mu_\tau(x)
=\int_{(0,b)}(t+\tau^2)^{N-1}\,d\nu(t).
$$
These prove necessity of the support and moment conditions.
This finite statement asserts neither an infinite orthogonal sequence
nor a functional on all polynomials.
\end{proof}

\begin{remark}\label{rem:single-cutoff-needed}
The simplicity proof uses the sign change rather than the sharp
indicator itself. It applies whenever each candidate $\lambda$ admits
$c_\lambda$ such that $(c_\lambda-t)(\chi(t)-\lambda)\,d\nu(t)$,
possibly with reversed sign, is positive with infinite support.
For input measures of infinite support, every bounded strictly monotone
$\chi:(0,b)\to(0,1)$ gives simple concentration spectrum for every $N$:
choose its sign-change separator $c_\lambda$ and repeat
\eqref{eq:positive-dual-cutoff} and the consecutive-zero argument.
Arbitrary indicators do not have this property. For $\nu=dt$ on $(0,2)$, $N=3$ and
$\chi=\mathbf1_{(1/2,3/2)}$, the bases
$\widehat R_0=1/\sqrt2$, $\widehat R_1=\sqrt{3/2}(t-1)$ and
$\widehat S_0=1/\sqrt2$ give
$\op{M}_\tau=\operatorname{diag}(1/2,1/2,1/8)$.
For fixed $N$, all moments are unnecessary: $d\nu(t)=t^{-(N+1)}\mathbf1_{(1,\infty)}(t)\,dt$ has a finite
moment of order $N-1$ and an infinite moment of order $N$.
\end{remark}

\begin{corollary}\label{cor:transfer}
For $N\ge2$, let $\op{K}_R,\op{K}_S$ be real symmetric matrices of orders
$N_R,N_S$, commuting with $\op{G}_R(\chi),\op{G}_S(\chi)$, respectively.
Then $\op{Q}^{\mathsf T}(\op{K}_R\oplus \op{K}_S)\op{Q}$ commutes with
every $\op{M}_\tau(\chi)$. Half-bandwidths $m_R,m_S$ become at most
$2\max(m_R,m_S)$. Every real symmetric commutant whose nonzero
diagonals have even offsets arises this way. If $\op{K}_R,\op{K}_S$ have simple spectra,
replacing $\op{K}_S$ by $\op{K}_S+c\Id$, where
$$
c=1+\|\op{K}_R\|_\infty+\|\op{K}_S\|_\infty,
\eqtag{eq:shift}
$$
makes the full spectrum simple, with distance at least one between
the two resulting block spectra. Here $\|\cdot\|_\infty$ is the
maximum absolute row sum.
\end{corollary}

\begin{proof}
Use \eqref{eq:isospectral}. Index differences double, and
$c-\|\op{K}_S\|_\infty>\|\op{K}_R\|_\infty$ separates the two spectra.
Conversely, the stated zero entries make
$\op{Q}\op{L}\op{Q}^{\mathsf T}$ block diagonal; its commutator
with $\op{G}_R\oplus \op{G}_S$ vanishes precisely when both scalar commutators do.
\end{proof}

\begin{proposition}\label{prop:normalized-transfer}
Under Corollary~\ref{cor:transfer}, suppose both scalar matrices have
simple spectra. For $j=R,S$, put
$$
\zeta_j\ge\max\{1,\|\op{K}_j\|_\infty\},\quad
\delta_j=\operatorname{gap}(\op{K}_j),
$$
where $\operatorname{gap}$ is the minimum spacing between consecutive
ordered eigenvalues, with value $\infty$ for a scalar matrix. Define
$$
\begin{aligned}
\widehat{\op{K}}_R&=\zeta_R^{-1}\op{K}_R-2\Id,\quad
\widehat{\op{K}}_S=\zeta_S^{-1}\op{K}_S+2\Id,\\[7pt]
\widehat{\op{L}}
&=\op{Q}^{\mathsf T}
\bigl(\widehat{\op{K}}_R\oplus\widehat{\op{K}}_S\bigr)\op{Q}.
\end{aligned}
\eqtag{eq:normalized-commutant}
$$
Then $\widehat{\op{L}}$ commutes with every $\op{M}_\tau(\chi)$,
has simple spectrum and satisfies the bandwidth bound of
Corollary~\ref{cor:transfer}. The distance between its two block spectra
is at least two, and
$$
\begin{aligned}
\|\widehat{\op{L}}\|_2&\le3,\\[7pt]
\operatorname{gap}(\widehat{\op{L}})
&\ge\min\{\delta_R/\zeta_R,\delta_S/\zeta_S,2\}.
\end{aligned}
\eqtag{eq:normalized-gap}
$$
Here $\|\op{A}\|_2=\max_{\|z\|_2=1}\|\op{A}z\|_2$.
For tridiagonal scalar matrices, $\widehat{\op{L}}$ has nonzero entries
only on the main and second diagonals.
\end{proposition}

\begin{proof}
The absolute row-sum bound gives
$|\lambda|\le\|\op{K}_j\|_\infty$ for each scalar eigenvalue.
Thus the eigenvalues of $\widehat{\op{K}}_R$ lie between $-3$ and $-1$,
and those of $\widehat{\op{K}}_S$ between $1$ and $3$.
The affine changes preserve commutation and divide the internal
spacings by $\zeta_j$. Orthogonal reconstruction preserves spectra
and the Euclidean operator norm. Corollary~\ref{cor:transfer} gives
the commutation and bandwidth assertions.
\end{proof}

This separation does not give a uniform lower bound for the internal
spacings as $N$ grows. Indeed,
$\operatorname{gap}(\widehat{\op{L}})\le6/(N-1)$ for $N\ge2$,
since its $N$ eigenvalues lie between $-3$ and $3$.
Certified upper bounds may be used for $\zeta_j$; powers of two permit
exact binary rescaling.

\begin{proposition}\label{prop:apriori-gap}
Under Proposition~\ref{prop:normalized-transfer}, suppose the scalar
matrices are tridiagonal with nonzero off-diagonal entries.
In the block $\op{K}_j$ of order $n_j\ge2$, write $d_{j,i}$ for the diagonal
entries ($1\le i\le n_j$) and $b_{j,i}$ for the links
($1\le i<n_j$), with $b_{j,0}=b_{j,n_j}=0$. Put
$$
\begin{aligned}
r_{j,i}&=|b_{j,i-1}|+|b_{j,i}|,\\[7pt]
W_j&=\max_i(d_{j,i}+r_{j,i})-\min_i(d_{j,i}-r_{j,i}),\\[7pt]
g_j&=\frac{2\prod_{i=1}^{n_j-1}|b_{j,i}|}
{\zeta_jW_j^{\,n_j-2}}.
\end{aligned}
\eqtag{eq:apriori-block-gap}
$$
For a scalar block set $g_j=\infty$. Then
$$
\operatorname{gap}(\widehat{\op{L}})\ge\min\{g_R,g_S,2\}>0.
\eqtag{eq:apriori-full-gap}
$$
If symmetric stored blocks $\op{B}_j$ satisfy
$\|\op{B}_j-\widehat{\op{K}}_j\|_2\le e_j$, the reconstructed matrix $\op{Q}^{\mathsf T}(\op{B}_R\oplus\op{B}_S)\op{Q}$
has gap at least
$$
\min\{g_R-2e_R,g_S-2e_S,2-e_R-e_S\},
\eqtag{eq:apriori-stored-gap}
$$
whenever this minimum is positive.
\end{proposition}

\begin{proof}
Let $\op{K}$ be one scalar block of order $n\ge2$, let
$p(\lambda)=\det(\lambda\Id-\op{K})$ and let $v_i$ be a unit
eigenvector at $\lambda_i$. The corner cofactor and the spectral
expansion of the adjugate give
$$
\prod_{k=1}^{n-1}|b_k|
=|p'(\lambda_i)v_i(1)v_i(n)|\le\tfrac12|p'(\lambda_i)|.
$$
If $\lambda_i$ is adjacent to a smallest gap $\delta$, then
$|p'(\lambda_i)|\le\delta W^{n-2}$ for any upper bound $W$ on the
spectral diameter. Testing an eigenvector at a coordinate of maximal
modulus places every eigenvalue between
$\min_i(d_i-r_i)$ and
$\max_i(d_i+r_i)$; hence the stated $W_j$ is such a bound.
Divide the resulting scalar gap by $\zeta_j$ and use the block
separation in Proposition~\ref{prop:normalized-transfer}.
A symmetric perturbation of norm $e_j$ moves each ordered eigenvalue
by at most $e_j$, which proves \eqref{eq:apriori-stored-gap}.
\end{proof}

These bounds use only matrix entries. They give a sufficient rounding
threshold before eigenvalue isolation; they need not be sharp for a
given input and may be very small as $N$ grows.

\section{Reconstructing commuting operators}

Take Laguerre or Jacobi input in Theorem~\ref{thm:deformation},
with $d\nu(t)=w(t)\,dt$ on $(0,b)$, where
$$
w(t)=
\begin{cases}
t^\beta e^{-t},\quad b=\infty,\quad\beta>-1,
&\text{Laguerre},\\[7pt]
t^\beta(1-t)^\alpha,\quad b=1,\quad\alpha,\beta>-1,
&\text{Jacobi}.
\end{cases}
\eqtag{eq:weights}
$$
Then $d\mu_\tau(x)=W_\tau(x)\,dx$ on $E_\tau$, with
$$
W_\tau(x)=|x+\tau|w(x^2-\tau^2).
\eqtag{eq:density}
$$
For $0<s<b$, let $\chi=\mathbf1_{(0,s)}$ be the indicator of $(0,s)$ and abbreviate
$\mathbf v_s=\mathbf v_\chi$ and $\mathbf u_{\tau,s}=\mathbf u_{\tau,\chi}$.
The coefficients in $(\sigma w)'=\eta w$ and the cutoff constants are
$$
\begin{array}{c|cc}
&\text{Laguerre}&\text{Jacobi}\\[7pt]
\sigma(t)&t&t(1-t)\\[7pt]
\eta(t)&\beta+1-t&\beta+1-(\alpha+\beta+2)t\\[7pt]
\kappa_R&r_R&r_R(r_R+\alpha+\beta+2)\\[7pt]
\kappa_S&r_S&r_S(r_S+\alpha+\beta+3).
\end{array}
$$

\begin{samepage}
\begin{proposition}\label{prop:commutation}
The operators
$$
\begin{aligned}
\op{T}_Rp
&=(t-s)\sigma p''+\bigl((t-s)\eta+\sigma\bigr)p'+\kappa_Rtp,\\[7pt]
\op{T}_Sp
&=(t-s)\sigma p''
+\bigl((t-s)(\eta+\sigma/t)+\sigma\bigr)p'+\kappa_Stp
\end{aligned}
\eqtag{eq:T}
$$
have symmetric tridiagonal matrices $\op{K}_R,\op{K}_S$ on
$\Pp_{r_R},\Pp_{r_S}$ in the $\widehat R,\widehat S$ bases, with positive
off-diagonal entries and hence simple spectra. Moreover,
$\op{K}_R\op{G}_R=\op{G}_R\op{K}_R$ and $\op{K}_S\op{G}_S=\op{G}_S\op{K}_S$.
\end{proposition}
\end{samepage}

\begin{proof}
Let $\partial p=p'$. For $\mathbf a=\mathbf v,\mathbf v_s$,
integration by parts gives
$$
-\partial'\bigl((x-s)\sigma\mathbf a\bigr)
=\bigl((x-s)\eta+\sigma\bigr)\mathbf a.
\eqtag{eq:dual-symmetry}
$$
At the full endpoints, $(t-s)\sigma(t)w(t)$ times any polynomial tends
to zero; at $s$, the cutoff factor vanishes.
For $f,g\in\Pp$, put $q=f'g-fg'$. Pairing
\eqref{eq:dual-symmetry} with $q$ gives
$$
\begin{aligned}
\langle\mathbf a,(\op{T}_Rf)g-f(\op{T}_Rg)\rangle
&=\langle\mathbf a,(x-s)\sigma q'
+\bigl((x-s)\eta+\sigma\bigr)q\rangle=0,\\[7pt]
\op{J}_{\mathbf a}\op{T}_R&=\op{T}_R'\op{J}_{\mathbf a}.
\end{aligned}
$$
For $\mathbf a=x\mathbf v,x\mathbf v_s$, replace $\eta$ by
$\eta+\sigma/x$; this proves both symmetries of $\op{T}_S$.
The coefficient of $t^{j+1}$ in $\op{T}_R(t^j)$ is $r_R-j$ for
Laguerre and $(r_R-j)(r_R+j+\alpha+\beta+2)$ for Jacobi.
For $\op{T}_S$, replace $r_R$ by $r_S$ and $\beta$ by $\beta+1$.
It is positive below the cutoff degree and zero at the cutoff.
Symmetry for the full functionals and this degree bound give symmetric tridiagonal
matrices with positive off-diagonal entries. In their orthonormal bases,
the truncated symmetries give
$$
\op{K}_R\op{G}_R=\op{G}_R\op{K}_R,\quad
\op{K}_S\op{G}_S=\op{G}_S\op{K}_S.
$$
An eigenvector of a tridiagonal matrix with nonzero off-diagonal entries
is determined by its first entry, so each symmetric block has simple spectrum.
\end{proof}

For $\op{A} p=\sigma p''+\eta p'$, $\op{X}p=tp$ and
$\eta=\eta_0+\eta_1t$, the operators in
\cite{Grunbaum,GVZ} can be written as
$$
\op{T}_R
=\tfrac12(\op{A}\op{X}+\op{X}\op{A})-s\op{A}
+\bigl(\kappa_R-\tfrac12\eta_1\bigr)\op{X}-\tfrac12\eta_0\Id,
\eqtag{eq:operator-product}
$$
with the corresponding substitutions for $\op{T}_S$.

\begin{corollary}\label{cor:local}
For the weights \eqref{eq:weights}, every $\tau\in\mathbb R$ and
$N\ge2$, the matrix $\op{M}_\tau(\mathbf1_{(0,s)})$ has a real symmetric
commuting matrix with only the main and second diagonals and simple
spectrum. It represents a second-order differential--reflection
operator preserving $\Pp_{N-1}$ and commuting on polynomials with
$$
(\op{C}_{\tau,N,s}F)(x)
=\sum_{j=0}^{N-1}\widehat P_j^{(\tau)}(x)
\langle\mathbf u_{\tau,s},\widehat P_j^{(\tau)}F\rangle.
\eqtag{eq:C}
$$
For $N=1$ the matrix is scalar.
On $\Pp_{N-1}$, \eqref{eq:C} agrees with the compressed
multiplier above and has matrix $\op{M}_\tau(\chi)$.
The normalized choice in Proposition~\ref{prop:normalized-transfer}
also represents an operator of order at most two with reflection.
\end{corollary}

\begin{proof}
For $F(x)=f(y)+(x-\tau)g(y)$, define
$$
(\op{T}_{\tau,N}F)(x)
=(\op{T}_Rf)(y)+(x-\tau)\bigl((\op{T}_S+c\Id)g\bigr)(y).
\eqtag{eq:lifted-operator}
$$
Its restriction to $\Pp_{N-1}$ has matrix
$\op{L}=\op{Q}^{\mathsf T}(\op{K}_R\oplus(\op{K}_S+c\Id))\op{Q}$,
with $c$ from \eqref{eq:shift}; Corollary~\ref{cor:transfer} gives simplicity.
The identities
$$
\begin{aligned}
g(y)&=\frac{F(x)-F(-x)}{2x},\\[7pt]
f(y)&=\frac{F(x)+F(-x)}2+\tau g(y)
\end{aligned}
$$
show that it is a second-order differential--reflection operator with ordinary leading
coefficient $(y-s)\sigma(y)/(4x^2)$.
The quotients for $f,g$ extend polynomially through $x=0$;
the action there is defined by \eqref{eq:lifted-operator}.
Let $\op{B}=\op{T}_R\oplus(\op{T}_S+c\Id)$.
For $\mathbf a=\mathbf u_\tau,\mathbf u_{\tau,s}$, the corresponding
dual pullback \eqref{eq:dual-pullback} or \eqref{eq:dual-mask} gives
$$
\begin{aligned}
\op{U}_\tau'\op{J}_{\mathbf a}\op{T}_{\tau,N}\op{U}_\tau
&=(\op{U}_\tau'\op{J}_{\mathbf a}\op{U}_\tau)\op{B}\\[7pt]
&=\op{B}'(\op{U}_\tau'\op{J}_{\mathbf a}\op{U}_\tau)
=\op{U}_\tau'\op{T}_{\tau,N}'\op{J}_{\mathbf a}\op{U}_\tau.
\end{aligned}
$$
Since $\op{U}_\tau$ and its transpose are invertible on the polynomial
spaces and their duals, $\op{T}_{\tau,N}$ is symmetric for both
functionals. In particular, for $j<N$ and every $F\in\Pp$,
$$
\langle\mathbf u_{\tau,s},\widehat P_j^{(\tau)}\op{T}_{\tau,N}F\rangle
=\langle\mathbf u_{\tau,s},
(\op{T}_{\tau,N}\widehat P_j^{(\tau)})F\rangle.
$$
Substitution in \eqref{eq:C}, using the cutoff preservation and the
symmetric matrix $\op{L}$, proves commutation on all polynomials.
The normalized realization is
$$
\widehat{\op{T}}_{\tau,N}
=\op{U}_\tau
\bigl((\zeta_R^{-1}\op{T}_R-2\Id)
\oplus(\zeta_S^{-1}\op{T}_S+2\Id)\bigr)\op{U}_\tau^{-1}.
$$
The formulas for $\op{U}_\tau^{-1}$ above show that its order is at
most two. Both dual symmetries and the cutoff preservation persist
under these affine changes, giving the same commutation argument.
\end{proof}
\begin{theorem}\label{thm:complete-pentadiagonal}
Let $N\ge2$ and let $\chi$ be a mask for which $\op{M}_\tau(\chi)$
has simple spectrum. Suppose each $\op{G}_j$, $j=R,S$, admits a real
symmetric tridiagonal
commutant $\op{K}_j$ with nonzero internal off-diagonal entries.
Then every real symmetric pentadiagonal matrix commuting with
$\op{M}_\tau(\chi)$ has the form
$$
\op{Q}^{\mathsf T}
\bigl((a_R\op{K}_R+b_R\Id)\oplus
      (a_S\op{K}_S+b_S\Id)\bigr)\op{Q},\quad
 a_R,b_R,a_S,b_S\in\mathbb R,
\eqtag{eq:all-pentadiagonal}
$$
Conversely every matrix in \eqref{eq:all-pentadiagonal} commutes
with $\op{M}_\tau$. This space has dimension $\min(N,4)$.
For $N\ge3$, half-bandwidth two is minimal among real symmetric
commuting matrices with simple spectrum in the degree order
$\widehat P_0^{(\tau)},\ldots,\widehat P_{N-1}^{(\tau)}$.
\end{theorem}

\begin{proof}
The simple full concentration spectrum gives simple spectra for
the two blocks and no eigenvalue in common. A commuting matrix therefore preserves each block: in
orthonormal eigenbases, an off-block entry would satisfy
$(\lambda_{R,i}-\lambda_{S,k})z_{ik}=0$.
A pentadiagonal matrix consequently has scalar tridiagonal blocks
$\op{A}_R,\op{A}_S$, with no entries on its first diagonals.

Each $\op{G}_j$ has simple spectrum, so $\op{A}_j$ and $\op{K}_j$
are diagonal in the same eigenbasis and commute with each other.
Write $b_i\ne0$ for the links of $\op{K}_j$ and $c_i$ for those of
$\op{A}_j$. Comparing second diagonals in their commutator gives
$b_ic_{i+1}=c_ib_{i+1}$. Thus $c_i=a_jb_i$ for a constant $a_j$.
The diagonal matrix $\op{A}_j-a_j\op{K}_j$ still commutes with
$\op{K}_j$; its adjacent diagonal entries must be equal.
Hence $\op{A}_j=a_j\op{K}_j+b_j\Id$.
This also holds for a scalar block, where the parameters are redundant.
The converse follows from Corollary~\ref{cor:transfer}.
Each block contributes two dimensions if its order is at least two,
and one otherwise, giving $\min(N,4)$.

A commuting matrix of half-bandwidth at most one must be diagonal,
by the same block preservation. Formula~\eqref{eq:all-pentadiagonal}
then makes it constant on each parity block. For $N\ge3$ at least
one block has order at least two, so its spectrum is not simple.
Corollary~\ref{cor:transfer} attains half-bandwidth two with simple
spectrum.
\end{proof}

For Laguerre and Jacobi input, both restricted supports are infinite
and Theorem~\ref{thm:concentration-multiplicities} gives simple spectrum.
Proposition~\ref{prop:commutation} supplies the scalar matrices.
Thus Theorem~\ref{thm:complete-pentadiagonal} applies for every interior cutoff.

\section{Applications}
\label{sec:CG}

We compute in the scalar coordinates of \eqref{eq:dual-mask}, using
$$
tp_k(t)=a_{k+1}p_{k+1}(t)+b_kp_k(t)+a_kp_{k-1}(t),\quad
p_{-1}=0,\quad a_0=0,
$$
with $p=\widehat R$ or $\widehat S$. These coefficients are evaluated
from the scalar input at the working precision; no monomial moment
matrix is factored. The scalar commutants
are formed from \eqref{eq:operator-product} or \eqref{eq:meixner-operator}.

For an $n\times n$ real symmetric tridiagonal matrix with diagonal $d_i$ and
off-diagonal $b_i$, define
$$
\begin{aligned}
D_0(\xi)&=1,\quad D_1(\xi)=d_1-\xi,\\[7pt]
D_k(\xi)&=(d_k-\xi)D_{k-1}(\xi)-b_{k-1}^2D_{k-2}(\xi),
\quad 2\le k\le n.
\end{aligned}
\eqtag{eq:sign-counts}
$$
Only finitely many endpoints are zeros of these polynomials.
At an endpoint where all $D_k(\xi)\ne0$, eliminating successive
coordinates gives a diagonal quadratic form with coefficients
$D_k(\xi)/D_{k-1}(\xi)$. The number of negative coefficients is the
number of eigenvalues below $\xi$. We call this the sign count.
Interval arithmetic verifies the signs; an undecidable endpoint is
replaced or evaluated at higher precision.

\refstepcounter{algorithm}
\label{alg:manufacture}
\begin{tcolorbox}[algobox,breakable,
  title={\textbf{Algorithm~\thealgorithm} Alternating concentration from scalar recurrences},
  title after break={\textbf{Algorithm~\thealgorithm} (continued)},
  colbacktitle=algogray,coltitle=black,fonttitle=\normalfont,
  titlerule=0.4pt]
\begin{algorithmic}[1]
\Require $N\ge1$, $\tau\in\mathbb R$, a cutoff $0<s<b$, scalar recurrence
formulas with initial constants $p_0$, symmetric tridiagonal commutants
with nonzero off-diagonal entries (scalar blocks allowed), a relative
value tolerance $\epsilon>0$ and, when eigenpolynomials are requested,
a tolerance $\epsilon_v>0$ for the sine of the angle.
\Ensure Enclosures of concentration values or deficits and, if requested,
unit eigenpolynomial approximations in the exact recurrence bases.
\State Set $N_R=\lceil N/2\rceil$, $N_S=\lfloor N/2\rfloor$ and
$y=x^2-\tau^2$. Omit empty blocks.
\For{each nonempty scalar block $j=R,S$}
\State Form the symmetric tridiagonal $\op{K}_j$ from the scalar recurrence.
\State Solve for approximate eigenpairs $(\eta_i,z_i)$ separately in this block.
\State For a scalar block, set $g_i=\infty$ and its angle bound to zero.
\For{each requested positive target $d_i$}
\Repeat
\State For a block of order greater than one, validate the non-target
distance $g_i$ for the current $\eta_i$ by sign counts.
\State Evaluate $p_{z_i}(t)=\sum_k(z_i)_kp_k(t)$ by recurrence.
Obtain $\widehat d_i$ from the requested positive integral
$$
q_{\mathrm{in}}(z_i)=
\frac{\int_{(0,s)}|p_{z_i}(t)|^2\,d\nu_j(t)}{\|z_i\|_2^2},\quad
q_{\mathrm{out}}(z_i)=
\frac{\int_{[s,b)}|p_{z_i}(t)|^2\,d\nu_j(t)}{\|z_i\|_2^2}.
$$
\State Bound the true normalized residual by $\rho_i$ and the integral
error by $e_{Q,i}$, including construction and rounding errors.
\State Set $B_i=e_{Q,i}+\Omega_j(\rho_i/g_i)^2$ as in
\eqref{eq:accuracy-bound}. In a one-dimensional block, omit the residual term.
\State If the tolerance fails, increase working precision, rebuild the
recurrence data and refine the eigenpair or integration enclosure.
\Until{$\widehat d_i>B_i$, $B_i\le\epsilon\widehat d_i/(1+\epsilon)$,
and $\rho_i/g_i\le\epsilon_v$ if eigenpolynomials are requested}
\EndFor
\EndFor
\State Enclose an extremal value by taking the maximum or minimum
of the respective interval endpoints over all candidates.
\State Reconstruct the unit approximations as
$p_{z_i}(y)/\|z_i\|_2$ in the $R$ block and
$(x-\tau)p_{z_i}(y)/\|z_i\|_2$ in the $S$ block; interlace the normalized coefficient vectors with
$\op{Q}^{\mathsf T}$.
\State \Return the scalar enclosures and reconstructed coefficient vectors.
\end{algorithmic}
\end{tcolorbox}
Here $d\nu_R=d\nu$ and $d\nu_S=t\,d\nu$.
The complement $[s,b)$ includes any mass at the cutoff.
Known zero values use absolute enclosures rather than a relative test.
In scalar blocks the angle test is automatic.
The scalar calculation is independent of $\tau$; \eqref{eq:norm}
transports its $L^2$ errors unchanged. For branch evaluation retain the
scalar variable $t$ and compute the smaller difference by
$r-|\tau|=t/(r+|\tau|)$, where $r=\sqrt{t+\tau^2}$.

The Laguerre case gives the $(-1)$-Meixner--Pollaczek family
\cite[Section~2.3]{Manufacture}, also described as deformed Hermite
polynomials in \cite[Section~7.2]{GVZChihara}; the Jacobi case gives the Chihara family
\cite[Proposition~5.1]{GVZChihara}.
The Chihara weight, branch cancellation and parameter-independent norms
are given in \cite[Section~5.1, (5.5)--(5.9)]{GVZChihara}.
For fixed $s$, all $\tau$ give the same concentration matrix on
$|\tau|<|x|<\sqrt{\tau^2+s}$. If instead the outer endpoint
$T=\sqrt{\tau^2+s}$ is fixed, the input cutoff is $s=T^2-\tau^2$.

Scalar Meixner concentration and its tridiagonal commutants were treated
in \cite{Perlstadt}. We use the operator construction in
\cite[Section~5]{GVZ}.

\begin{example}\label{ex:meixner}
Let $\beta>0$, $0<q<1$ and let $m\ge0$ be an integer.
On the grid $\{\beta+k:k=0,1,\ldots\}$, set
$$
\nu_j=\sum_{k=0}^{\infty}\frac{(\beta+j)_kq^k}{k!}\,
\delta_{\beta+k},\quad j=0,1.
\eqtag{eq:meixner-measures}
$$
Here $\delta_a$ is the unit point mass at $a$,
$(a)_k=a(a+1)\cdots(a+k-1)$ and $(a)_0=1$.
Both measures have all moments finite, and
$t\,d\nu_0(t)=\beta\,d\nu_1(t)$.
\begin{samepage}
The ordinary Meixner operators \cite[Section~18.22(ii)]{DLMF} are
$$
\op{A}_jp(t)
=q(t+j)\bigl(p(t+1)-p(t)\bigr)
 +(t-\beta)\bigl(p(t-1)-p(t)\bigr).
$$
\end{samepage}
Take $b=\infty$, $\nu=\nu_0$, $\chi=\mathbf1_{(0,s)}$,
$s=\beta+m+\tfrac12$, $N\ge2$ and $r_0=r_R$, $r_1=r_S$.
With $\op{X}p=tp$, put
$$
\op{H}_j
=\tfrac12(\op{A}_j\op{X}+\op{X}\op{A}_j)
-s\op{A}_j+(1-q)(r_j+\tfrac12)\op{X}.
\eqtag{eq:meixner-operator}
$$
Its coefficient raising degree $n$ is $(1-q)(r_j-n)$.
Writing $w_j(k)=(\beta+j)_kq^k/k!$, the balance identity
$w_j(k)q(k+\beta+j)=(k+1)w_j(k+1)$ cancels adjacent summands.
The coefficient $q(t+j)(t+\tfrac12-s)$ of $p(t+1)$ vanishes at
$t=\beta+m$, and the lower link vanishes at $t=\beta$.
For nonzero polynomials $f,g$, the boundary term at infinity is
$O(q^k k^{\beta+j+1+\deg f+\deg g})$ and tends to zero.
Consequently the full and truncated sums give
$$
\op{J}_{\mathbf a}\op{H}_j=\op{H}_j'\op{J}_{\mathbf a},\quad
\mathbf a=x^j\mathbf v,x^j\mathbf v_s,\quad j=0,1.
$$
Let $\op{K}_R,\op{K}_S$ be the matrices of $\op{H}_0,\op{H}_1$ in the
$\widehat R,\widehat S$ bases. They are symmetric tridiagonal with
positive off-diagonal entries, have simple spectra, and commute with
$\op{G}_R,\op{G}_S$, respectively.

The reconstructed measure is supported at
$$
x=\pm\sqrt{\beta+k+\tau^2},\quad k=0,1,\ldots.
$$
The mask retains $k=0,\ldots,m$. The blocks $\op{G}_R,\op{G}_S$ have ranks
$\min(N_R,m+1),\min(N_S,m+1)$, respectively.
The full concentration spectrum is simple precisely when
$m+1\ge\lfloor N/2\rfloor$. Under this condition,
Theorem~\ref{thm:complete-pentadiagonal} classifies all its real
symmetric pentadiagonal commutants.
The measure and orthonormal polynomials
are given by \eqref{eq:measure} and \eqref{eq:deformed-polynomials}.
For $\tau=0$, this is the Meixner system of
\cite[Section~3, Case~3]{Jafarov}, up to normalization and a scaling
of the variable.
Corollary~\ref{cor:transfer} gives a commuting matrix with only the main
and second diagonals and simple spectrum after \eqref{eq:shift}.

Write $\op{M}=\op{M}_\tau(\chi)$, which is independent of $\tau$.
For $\beta=3/2$, $m=5$ and $N=12$, the extremal concentration values are
$$
\begin{array}{ccc}
q&\lambda_{\min}(\op{M})&1-\lambda_{\max}(\op{M})\\[7pt]
1/100&0.636574&6.43713\times10^{-27}\\[7pt]
1/3&5.32886\times10^{-6}&2.79643\times10^{-9}\\[7pt]
99/100&2.86392\times10^{-31}&0.884905.
\end{array}
$$
The values were computed at 100 digits and checked at 130.

The relative commutator norms
$\|\op{K}_j\op{G}_j-\op{G}_j\op{K}_j\|_F/(\|\op{K}_j\|_F\|\op{G}_j\|_F)$ were below
$3\times10^{-101}$ at 100 digits, for $j=R,S$.
Here $\|\op{A}\|_F^2=\sum_{i,k}|\op{A}_{ik}|^2$.
\end{example}

\begin{remark}\label{rem:rank-deficient-commutant}
The simplicity condition in Theorem~\ref{thm:complete-pentadiagonal}
cannot be omitted in general. In Example~\ref{ex:meixner}, take
$\beta=1$, $q=1/2$, $m=0$ and $N=5$. Then
$$
\op{G}_R=rr^{\mathsf T},\quad
r=(1/\sqrt2,-1/2,1/(2\sqrt2))^{\mathsf T},\quad
\op{A}=\begin{pmatrix}1/\sqrt2&1&0\\[7pt]1&\sqrt2&0\\[7pt]0&0&0\end{pmatrix}
$$
satisfies $\op{A}r=0$. Thus $\op{A}$ commutes with $\op{G}_R$,
but it is not affine in $\op{K}_R$: its second link is zero and its
first is nonzero, whereas both links of $\op{K}_R$ are nonzero.
The matrix $\op{K}_R+\op{A}$ remains irreducible tridiagonal.
A sufficiently large shift of $\op{K}_S$ therefore gives a
simple-spectrum pentadiagonal commutant outside
\eqref{eq:all-pentadiagonal}.
\end{remark}

\begin{remark}\label{rem:CG}
Let $\mu>-1/2$, with $T>0$ for Laguerre and
$\alpha>-1$, $0<T<1$ for Jacobi.
Set $\tau=0$, $s=T^2$, $\beta=\mu-\tfrac12$ and
$\chi=\mathbf1_{(0,s)}$. Apply the alternating reconstruction of
Theorem~\ref{thm:deformation} to the scalar operators of
Proposition~\ref{prop:commutation}. The weights \eqref{eq:weights} give
$$
(\op{L}_R,\op{L}_S)=
\begin{cases}
\bigl(2\op{K}_R,\;2\op{K}_S+(s+1)\Id\bigr),&\text{Laguerre},\\[7pt]
\bigl(4\op{K}_R,\;4\op{K}_S+2\bigl(\mu+1+(\alpha+\mu+1)s\bigr)\Id\bigr),
&\text{Jacobi}.
\end{cases}
\eqtag{eq:exact}
$$
The matrix $\op{L}=\op{Q}^{\mathsf T}(\op{L}_R\oplus \op{L}_S)\op{Q}$ is
\cite[(6), (14)]{CG} in the positive-leading-coefficient orthonormal basis.
For $N\ge2$, these are particular choices in the complete family
\eqref{eq:all-pentadiagonal}; the separated choice
\eqref{eq:normalized-commutant} has simple spectrum for every admissible
parameter value.
For $N\ge2$, let $\op{D}_{\mathrm H},\op{D}_{\mathrm G}$ denote the
operators \cite[(7), (15)]{CG}, and write $\op{T}_\mu F(x)=F'(x)+\mu(F(x)-F(-x))/x$.
For $t=x^2$, one has $\op{T}_\mu f(t)=2xf'(t)$ and
$\op{T}_\mu(xg(t))=2tg'(t)+(2\mu+1)g(t)$.
With $N_0=2r_R$ and $N_1=2r_S+1$, substitution gives, on
$\Pp\oplus\Pp$,
$$
\begin{aligned}
\op{U}_0^{-1}\op{D}_{\mathrm H}\op{U}_0
&=2\op{T}_R\oplus\bigl(2\op{T}_S+(s+1)\Id\bigr),\\[7pt]
\op{U}_0^{-1}\op{D}_{\mathrm G}\op{U}_0
&=4\op{T}_R\oplus\bigl(4\op{T}_S
+2\bigl(\mu+1+(\alpha+\mu+1)s\bigr)\Id\bigr).
\end{aligned}
$$
The first row uses the Laguerre operators and the second the Jacobi
operators. Taking matrices with the recurrences \cite[(4), (11)]{CG}
gives the entries in \cite[(6), (14)]{CG}. These identities give the
polynomial actions asserted in \cite[Theorems~2.1 and~3.1]{CG},
with the basis conventions below. Their commutation
follows from \eqref{eq:dual-symmetry},
the alternating pullback \eqref{eq:dual-mask} and the cutoff preservation,
as in Corollary~\ref{cor:local}.
For Gegenbauer, put $c=\alpha+\mu+\tfrac12$. When $-1<c<0$, the
polynomials normalized in \cite{CG} have negative leading coefficients
in positive degrees. Our formulas use positive-leading orthonormal bases.
Passing to their literal basis conjugates both $\op{L}$ and $\op{M}$
by $\op{J}=\operatorname{diag}(1,-1,\ldots,-1)$.
At $c=0$ we use the normalized Jacobi limit. In the recurrences
\cite[(4), (11)]{CG}, we set $a_0=0$ and evaluate removable quotients
by continuity.
For $N=1$ the matrix identity is scalar.
\end{remark}

The Hermite blocks in \eqref{eq:exact} also determine every cutoff
at which the unmodified full matrix has repeated eigenvalues.
For $N\ge2$, put $a=\mu+\tfrac12>0$, $s=T^2$ and
$m_0=N_R$, $m_1=N_S$. Write $\op{L}_0=\op{L}_R$ and
$\op{L}_1=\op{L}_S$. In parity $\varepsilon=0,1$, their diagonal
entries $d_{\varepsilon,j}$ and positive off-diagonal entries
$b_{\varepsilon,j}$ satisfy
$$
\begin{aligned}
d_{\varepsilon,j}
&=2(m_\varepsilon-1-j)(2j+\varepsilon+a)
 +(\varepsilon+2j)(s+1),\\[7pt]
b_{\varepsilon,j}^2
&=4(m_\varepsilon-1-j)^2(j+1)(j+a+\varepsilon).
\end{aligned}
\eqtag{eq:Hermite-block-entries}
$$
Here $0\le j<m_\varepsilon$ in the first row and
$0\le j<m_\varepsilon-1$ in the second. Define the monic polynomials
$$
p_R(\lambda;a,s)=\det(\lambda\Id-\op{L}_R),\quad
p_S(\lambda;a,s)=\det(\lambda\Id-\op{L}_S).
$$

\begin{theorem}\label{thm:Hermite-exceptions}
For $N\ge2$, the resultant
$$
\Delta_N(a,s)=\operatorname{Res}_\lambda(p_R,p_S)
\eqtag{eq:exceptional-polynomial}
$$
is a polynomial in $a,s$ with integer coefficients. Its degree in $s$
is $m_0m_1=\lfloor N^2/4\rfloor$, with leading coefficient
$$
c_N=\prod_{i=0}^{m_0-1}\prod_{j=0}^{m_1-1}(2i-2j-1)\ne0.
$$
The full Hermite matrix $\op{L}$ has simple spectrum precisely when
$\Delta_N(a,s)\ne0$. Thus, for each fixed $\mu>-1/2$, at most
$\lfloor N^2/4\rfloor$ positive cutoffs $s=T^2$ have repeated
eigenvalues. In particular,
$$
s>2\|\op{L}(a,0)\|_\infty
\eqtag{eq:large-cutoff}
$$
guarantees simple spectrum.
\end{theorem}

\begin{proof}
The tridiagonal determinant recurrence is
$$
\begin{aligned}
P_{\varepsilon,0}&=1,\quad
P_{\varepsilon,1}=\lambda-d_{\varepsilon,0},\\[7pt]
P_{\varepsilon,j+1}
&=(\lambda-d_{\varepsilon,j})P_{\varepsilon,j}
 -b_{\varepsilon,j-1}^2P_{\varepsilon,j-1}
 \quad(1\le j<m_\varepsilon).
\end{aligned}
$$
It gives $p_R=P_{0,m_0}$ and $p_S=P_{1,m_1}$ with integer
coefficients in $a,s,\lambda$. The resultant is a determinant in these
coefficients, so its coefficients are integers too. Equivalently,
if $\lambda_{\varepsilon,j}$ are the increasingly ordered
eigenvalues of $\op{L}_\varepsilon$,
$$
\Delta_N(a,s)
=\prod_{i=0}^{m_0-1}\prod_{j=0}^{m_1-1}
  (\lambda_{0,i}-\lambda_{1,j}).
$$
Write $\op{L}_\varepsilon=\op{B}_\varepsilon+s\op{D}_\varepsilon$,
where $\op{B}_\varepsilon=\op{L}_\varepsilon(a,0)$ and
$$
\op{D}_\varepsilon
=\operatorname{diag}(\varepsilon,\varepsilon+2,\ldots,
                     \varepsilon+2m_\varepsilon-2).
$$
For each fixed $a$,
$|\lambda_{\varepsilon,j}-(\varepsilon+2j)s|
\le\|\op{B}_\varepsilon\|_2$, so
$\lambda_{\varepsilon,j}=(\varepsilon+2j)s+O(1)$ as $s\to\infty$.
Consequently,
$$
\Delta_N(a,s)=c_Ns^{m_0m_1}+O(s^{m_0m_1-1}).
$$
Since this holds for every $a>0$, it proves the polynomial degree
and the constant leading coefficient. The positive off-diagonal
entries give simple spectra within the blocks. Thus a repeated
eigenvalue occurs precisely when $p_R,p_S$ have a common root.
The degree bound gives the finite exceptional set.
Finally, $\op{L}(a,s)=\op{L}(a,0)+s\operatorname{diag}(0,1,\ldots,N-1)$.
Its ordered eigenvalues differ from $0,s,\ldots,(N-1)s$ by at most
$\|\op{L}(a,0)\|_2\le\|\op{L}(a,0)\|_\infty$, proving
\eqref{eq:large-cutoff}.
\end{proof}

For $N=3$, the recurrence gives
$$
\Delta_3(a,s)=2a(s-1)-(s+1)^2.
$$
For $a>4$ its roots $s=a-1\pm\sqrt{a(a-4)}$ are positive and
distinct, so the exceptional-cutoff bound is attained.
For $a=4$ the sole exceptional cutoff is $s=3$;
for $0<a<4$ there are none.

\begin{proposition}\label{prop:odd-Hermite-collisions}
For every odd $N\ge3$ and every $S>0$, there are parameters $a>0$
and $s>S$ for which the full Hermite matrix has an eigenvalue of
multiplicity two.
\end{proposition}

\begin{proof}
Write $N=2m+1$ and let $\ell_R(a,s)$ and $\ell_S(a,s)$ be the
smallest eigenvalues of the two blocks. Fix $s>1$ and isolate the
last coordinate of either block as
$$
\op{H}_a=
\begin{pmatrix}
\op{A}_a&b_a\\[7pt]
b_a^{\mathsf T}&d
\end{pmatrix}.
$$
For $m=1$ the odd block is the scalar $s+1$ and needs no decomposition.
Formula \eqref{eq:Hermite-block-entries} gives
$$
\begin{aligned}
\frac{\op{A}_{R,a}}{2a}
&\longrightarrow\operatorname{diag}(m,\ldots,1),\quad
\frac{b_{R,a}}{\sqrt{2a}}\longrightarrow\sqrt{2m}e_m,\\[7pt]
\frac{\op{A}_{S,a}}{2a}
&\longrightarrow\operatorname{diag}(m-1,\ldots,1)\quad(m\ge2),\\[7pt]
\frac{b_{S,a}}{\sqrt{2a}}
&\longrightarrow\sqrt{2(m-1)}e_{m-1}\quad(m\ge2).
\end{aligned}
$$
Here $a\to\infty$, $e_j$ is the last coordinate vector in
$\mathbb R^j$, $d_R=2m(s+1)$ and $d_S=(2m-1)(s+1)$.
The smallest eigenvalue of each nonempty $\op{A}_a$ tends to infinity.
For sufficiently large $a$, completing the square gives
$$
\min\{0,d-b_a^{\mathsf T}\op{A}_a^{-1}b_a\}
 \le\lambda_{\min}(\op{H}_a)\le d.
$$
The lower bounds stay finite, so these eigenvalues remain bounded.
Eliminating the first block gives
$$
\lambda=d-b_a^{\mathsf T}(\op{A}_a-\lambda\Id)^{-1}b_a.
$$
The limits above yield
$$
\ell_R(a,s)\longrightarrow2ms,\quad
\ell_S(a,s)\longrightarrow(2m-1)s+1.
\eqtag{eq:odd-block-limits}
$$
Choose $s_0>\max\{1,S\}$. For sufficiently large finite $a$,
$\ell_R(a,s_0)-\ell_S(a,s_0)>0$. Hold this $a$ fixed.
Testing the first coordinate in the even block and using
$\op{D}_1\ge\Id$ in the odd block gives
$$
\ell_R(a,s)\le2ma,\quad
\ell_S(a,s)\ge s+\ell_S(a,0).
$$
Thus the difference is negative for sufficiently large $s$.
The difference therefore vanishes at some $s>s_0$.
Each block has simple spectrum, so their common smallest eigenvalue
has multiplicity two in the full matrix.
\end{proof}

The parameters in this proposition depend on $N$ and $S$.
Proposition~\ref{prop:normalized-transfer} separates the Hermite blocks
for every admissible $N,\mu,T$.

For the Jacobi reconstruction, the exceptional cutoffs also form a finite
set, including on the parameter line where two diagonal slopes coincide. Fix
$\alpha>-1$, $\mu>-1/2$, put $\beta=\mu-\tfrac12$, $s=T^2$ and
$$
\vartheta=2\alpha+2\mu+1>-2,\quad
m_0=N_R,\quad m_1=N_S,\quad K=m_0m_1=\lfloor N^2/4\rfloor,
$$
and use the Jacobi blocks of \eqref{eq:exact}. Write $u_n>0$ ($n\ge1$) for the full recurrence coefficients, so that
$$
x\widehat P_n^{(0)}
=u_{n+1}\widehat P_{n+1}^{(0)}+u_n\widehat P_{n-1}^{(0)},\quad
u_0=0,\quad\widehat P_{-1}^{(0)}=0.
$$
Define
$$
p_R(\lambda;s)=\det(\lambda\Id-\op{L}_R(s)),\quad
p_S(\lambda;s)=\det(\lambda\Id-\op{L}_S(s)),\quad
\Xi_N(s)=\operatorname{Res}_\lambda(p_R,p_S).
$$

\begin{theorem}\label{thm:Jacobi-exceptions}
For every $N\ge2$ and every fixed admissible pair $\alpha,\mu$,
$\Xi_N$ is a nonzero real polynomial. The full Jacobi matrix of
\cite[(14)]{CG} has simple spectrum precisely when $\Xi_N(T^2)\ne0$.
Its degree is $K$ if $\vartheta\ne-1$. If $\vartheta=-1$, its degree is $K-1$,
except that for
$$
N=2m+1\ge3,\quad
\mu=-\frac{6m}{8m^2+8m-1}
$$
the degree is $K-2$. Thus the number of exceptional cutoffs
$0<T<1$ is at most the corresponding degree. Every repeated
eigenvalue has multiplicity two.
\end{theorem}

\begin{proof}
Write $\op{L}(s)=\op{A}+s\op{D}$, where
$$
\op{A}=\op{L}(0),\quad
\op{D}=\operatorname{diag}(q_0,\ldots,q_{N-1}),\quad
q_n=n(n+\vartheta).
$$
The entries on the second diagonals are independent of $s$ and
strictly positive. Indeed, in \cite[(14)]{CG} they are
$\Lambda_n u_{n+1}u_{n+2}$, with
$$
k_n=\left\lfloor\frac{N-1-n}{2}\right\rfloor,\quad
\Lambda_n=4k_n\left(n+k_n+\frac{\vartheta+2}{2}\right)
\quad(0\le n<N).
$$
For the internal links $0\le n\le N-3$, one has $\Lambda_n>0$.
Here the recurrence coefficients are taken in the
positive-leading orthonormal normalization, with removable quotients
evaluated by continuity. Hence each parity block is an irreducible
symmetric tridiagonal matrix and has simple spectrum. Its determinant
recurrence gives a polynomial in $s,\lambda$. The resultant therefore
is a polynomial in $s$, and vanishes exactly when the blocks have a
common eigenvalue. Such an eigenvalue has multiplicity two in the
full matrix.

The slopes are strictly increasing within each parity. The eigenvalue
in parity $\varepsilon$ corresponding to the slope $q_{2j+\varepsilon}$
is $q_{2j+\varepsilon}s+O(1)$ as $s\to\infty$. Consequently, if
$\vartheta\ne-1$, the leading coefficient of the resultant is
$$
\prod_{i=0}^{m_0-1}\prod_{j=0}^{m_1-1}
 (q_{2i}-q_{2j+1})
=\prod_{i=0}^{m_0-1}\prod_{j=0}^{m_1-1}
 (2i-2j-1)(2i+2j+1+\vartheta)\ne0.
$$
This proves that its degree is $K$.

Suppose now that $\vartheta=-1$. Then $-1/2<\mu<0$, and the sole coincidence
between slopes of different parities is $q_0=q_1=0$. Put
$A_n=\op{A}_{n,n}$. The recurrence coefficients needed below are
$$
u_1^2=1+2\mu,\quad u_2^2=-\frac{4\mu}{3},\quad
u_3^2=\frac{3+2\mu}{15}.
$$
For $N=2m$, direct substitution in \cite[(14)]{CG} gives
$$
A_1-A_0=4m-2+2\mu-\frac{16\mu}{3}(m-1)^2>0.
$$
All other cross-parity eigenvalue differences grow linearly with
$s$. The exceptional pair tends to $A_0-A_1\ne0$, so the resultant
has degree $K-1$.

For $N=2m+1$, substitution instead gives
$$
A_1-A_0=-4m-\frac{2\mu}{3}(8m^2+8m-1).
$$
It is nonzero except at the parameter stated in the theorem; this
again gives degree $K-1$. At that parameter, denote by $\eta_0(s)$
and $\eta_1(s)$ the eigenvalues belonging to the zero slopes.
Isolating the first coordinate and eliminating the remaining coordinates gives
$$
\eta_0(s)=A_0-\frac{C_0}{s}+O(s^{-2}),\quad
\eta_1(s)=A_1-\frac{C_1}{s}+O(s^{-2}),
$$
where
$$
\begin{aligned}
C_0&=\frac{\Lambda_0^2u_1^2u_2^2}{q_2}
 =-\frac{8\mu}{3}m^2(2m+1)^2(1+2\mu),\\[7pt]
C_1&=\frac{\Lambda_1^2u_2^2u_3^2}{q_3}
 =-\frac{8\mu}{135}(m-1)^2(2m+1)^2(3+2\mu).
\end{aligned}
$$
For $m=1$ the odd block is scalar and $C_1=0$. The same formula
therefore covers all $m\ge1$. At the specified parameter one has
$-2/5\le\mu<0$, since
$$
8m^2+8m-1-15m=(m-1)(8m+1)\ge0.
$$
It follows that
$$
45m^2(1+2\mu)-(m-1)^2(3+2\mu)
\ge 9m^2-3(m-1)^2>0,
$$
so $C_0>C_1$. Thus the exceptional eigenvalue difference is a
nonzero multiple of $s^{-1}$, whereas all the other $K-1$
cross-parity differences are nonzero multiples of $s$ to leading
order. The resultant has degree exactly $K-2$. This also proves
that it never vanishes identically.
\end{proof}

For example, at $N=3$, $\alpha=-3/5$ and $\mu=-2/5$, the resultant
is the nonzero constant $\Xi_3(s)=-96/25$. The reduction of its
degree therefore creates no persistent degeneracy.

The scalar residual bounds are as follows \cite{DavisKahan}.
\begin{proposition}\label{prop:numerical}
Under Corollary~\ref{cor:transfer}, suppose each $\op{K}_j$, $j=R,S$, has
simple spectrum, and let $\op{V}_j$ be an orthogonal matrix of its eigenvectors.
Then $\op{D}_j=\op{V}_j^{\mathsf T}\op{G}_j\op{V}_j$ is diagonal and, for every $\tau$,
$$
\op{Y}=\op{Q}^{\mathsf T}(\op{V}_R\oplus \op{V}_S),\quad
\op{M}_\tau(\chi)\op{Y}=\op{Y}(\op{D}_R\oplus \op{D}_S).
\eqtag{eq:numerical-transfer}
$$
The spectra of $\op{K}_R$ and $\op{K}_S$ may overlap.

\begin{samepage}
In a block of size at least two, let $(\kappa_i,v_i)$ be an eigenpair
of $\op{K}_j$, with $\|v_i\|_2=1$, and set
$\lambda_i=v_i^{\mathsf T}\op{G}_jv_i$. For a real unit vector $z$ and $\eta\in\mathbb R$,
put $r=\op{K}_jz-\eta z$ and
$\delta=\min_{k\ne i}|\kappa_k-\eta|$.
Write $\omega_j=\lambda_{\max}(\op{G}_j)-\lambda_{\min}(\op{G}_j)$;
then $\omega_j\le\min\{1,\operatorname{tr}\op{G}_j,
\operatorname{tr}(\Id-\op{G}_j)\}$.
If $\delta>0$, then
$$
\sin\angle(z,v_i)\le\frac{\|r\|_2}{\delta},\quad
|z^{\mathsf T}\op{G}_jz-\lambda_i|
\le\omega_j\frac{\|r\|_2^2}{\delta^2}.
\eqtag{eq:residual-bound}
$$
Both bounds are preserved by reconstruction with $\op{U}_\tau$.
Here $\angle$ denotes the acute angle between the spans of the vectors.

\end{samepage}
For a nonzero stored vector $z$, let $d_i$ be either $\lambda_i$ or
$1-\lambda_i$. Suppose $\rho$ bounds the true residual
$\|\op{K}_jz-\eta z\|_2/\|z\|_2$, $0<g\le\delta$,
$\Omega_j\ge\omega_j$, and $\widehat d$ approximates the corresponding
normalized positive form with absolute error at most $e_Q$.
Then
$$
\begin{aligned}
|\widehat d-d_i|&\le B=e_Q+\Omega_j(\rho/g)^2,\\[7pt]
\left|\frac{\widehat d}{d_i}-1\right|
&\le\frac{B}{\widehat d-B}\quad(\widehat d>B).
\end{aligned}
\eqtag{eq:accuracy-bound}
$$
For a positive target, $B\le\epsilon\widehat d/(1+\epsilon)$ with
$\widehat d>B$ guarantees relative error
at most $\epsilon$. Construction and rounding errors must be included
in $\rho$ and $e_Q$; $\Omega_j=1$ is always admissible.
\end{proposition}

\begin{proof}
Commutation and simplicity make every eigenvector of $\op{K}_j$ an
eigenvector of $\op{G}_j$. This proves \eqref{eq:numerical-transfer}.
Expand $z=\sum_k a_kv_k$ in their common orthonormal eigenbasis.
Since $0\le \op{G}_j\le\Id$,
$$
\begin{aligned}
\|r\|_2^2
&=\sum_k|\kappa_k-\eta|^2|a_k|^2
\ge\delta^2\sum_{k\ne i}|a_k|^2,\\[7pt]
|z^{\mathsf T}\op{G}_jz-\lambda_i|
&=\left|\sum_{k\ne i}(\lambda_k-\lambda_i)|a_k|^2\right|
\le\omega_j\sum_{k\ne i}|a_k|^2.
\end{aligned}
$$
Use \eqref{eq:norm} for the lifted vectors. In dimension one, both errors vanish.
Apply the same argument to $z/\|z\|_2$ and to $\Id-\op{G}_j$ for a
deficit, then add $e_Q$. Finally $d_i\ge\widehat d-B>0$ gives the
relative bound. If $\|\op{K}_j-\widetilde{\op{K}}_j\|_2\le e_K$, one may use
$\rho\ge\|\widetilde{\op{K}}_jz-\eta z\|_2/\|z\|_2+e_K$.
\end{proof}

The affine changes in Proposition~\ref{prop:normalized-transfer}
preserve the scalar eigenvectors. Applied also to the approximate
eigenvalue, they divide both the residual and the within-block spectral
distance by $\zeta_j$, leaving their ratio unchanged.
\begin{proposition}\label{prop:rounding-certificate}
Under Proposition~\ref{prop:normalized-transfer}, let $\op{B}_j$ be
stored real symmetric approximations to $\widehat{\op{K}}_j$,
with the same band structure, satisfying
$$
\|\op{B}_j-\widehat{\op{K}}_j\|_2\le e_j,\quad e_j\ge0,\quad j=R,S.
$$
Suppose a complete ordered list of closed intervals
$I_{j,k}$, with endpoints $a_{j,k},b_{j,k}$, each containing exactly
one eigenvalue of $\op{B}_j$, has been certified. For tridiagonal
$\op{B}_j$, use \eqref{eq:sign-counts}. Put
$$
\gamma_j=\min_{1\le k<N_j}(a_{j,k+1}-b_{j,k}),
$$
with $\gamma_j=\infty$ for a scalar block. If
$$
\gamma_j>2e_j\quad(j=R,S),\quad e_R+e_S<2,
$$
then $\widetilde{\op{L}}
=\op{Q}^{\mathsf T}(\op{B}_R\oplus\op{B}_S)\op{Q}$ has simple spectrum,
and
$$
\begin{aligned}
\operatorname{gap}(\widetilde{\op{L}})
&\ge\min\{\gamma_R,\gamma_S,2-e_R-e_S\},\\[7pt]
\operatorname{gap}(\widehat{\op{L}})
&\ge\min\{\gamma_R-2e_R,\gamma_S-2e_S,2\}.
\end{aligned}
\eqtag{eq:stored-gaps}
$$
In a block of order at least two, for a computed pair $(\xi,z)$
with $\xi\in I_{j,i}$ and $z\ne0$, take certified bounds
$$
\begin{aligned}
\rho&\ge\frac{\|\op{B}_jz-\xi z\|_2}{\|z\|_2}+e_j,\\[7pt]
0<g&\le\min_{k\ne i}\operatorname{dist}(\xi,I_{j,k})-e_j.
\end{aligned}
\eqtag{eq:stored-residual}
$$
The angle and positive-form bounds of Proposition~\ref{prop:numerical}
then apply with these $\rho,g$. For a scalar block the angle is zero
and the positive-form bound uses $B=e_Q$. Construction, normalization and
rounding errors are included in $e_j$; residual-evaluation errors are
included in $\rho$.
The stored matrix need not commute exactly with $\op{M}_\tau(\chi)$,
but, with $\omega_j$ as in Proposition~\ref{prop:numerical},
$$
\|\widetilde{\op{L}}\op{M}_\tau-\op{M}_\tau\widetilde{\op{L}}\|_2
\le\max_{j=R,S}(\omega_je_j)\le\max\{e_R,e_S\}.
\eqtag{eq:stored-commutator}
$$
\end{proposition}

\begin{proof}
The ordered eigenvalues of $\widehat{\op{K}}_j$ differ from those
of $\op{B}_j$ by at most $e_j$. Their intervals expanded by $e_j$
therefore remain disjoint, proving the exact-block gap bounds.
The stored $R$ spectrum lies between $-3-e_R$ and $-1+e_R$,
and the stored $S$ spectrum between $1-e_S$ and $3+e_S$.
This proves \eqref{eq:stored-gaps}. The residual for the exact block
is bounded by $\rho$, and its distance to every non-target eigenvalue
is bounded below by $g$, proving \eqref{eq:stored-residual}.
For the commutator, use the alternating decomposition
\eqref{eq:isospectral}. In each scalar block, center $\op{G}_j$
at the midpoint of its spectral extremes. Its centered norm is
$\omega_j/2$, so the commutator norm is at most $\omega_je_j$.
\end{proof}

For tridiagonal data, a certified upper bound for $e_j$ is the largest
row sum of the three entry-error bounds. All quantities in this
certificate can therefore be obtained from recurrence enclosures and
sign counts. If an exact eigenvector of the stored block is used,
the sine of its angle from the corresponding exact eigenvector is at most
$e_j/(\gamma_j-e_j)$.

\begin{corollary}\label{cor:finite-termination}
Under the input assumptions of Algorithm~\ref{alg:manufacture}, fix a
finite $N$. Suppose the initial constants, recurrence coefficients,
scalar commutant entries and positive-form evaluations at computable
coefficient vectors admit computable enclosures of arbitrarily small
width. With scalar eigenpairs refined by sign counts and
the tridiagonal recurrence, Algorithm~\ref{alg:manufacture} terminates
for every strictly positive target, every $\epsilon>0$ and every
requested $\epsilon_v>0$. After a sign choice, the reconstructed
unit eigenpolynomial, represented in the exact orthonormal bases, has
$L^2$ distance at most $\sqrt2\epsilon_v$ from its exact counterpart,
for every $\tau$.
\end{corollary}

\begin{proof}
Each scalar eigenvalue is simple and has positive distance from its
neighbors. Sign counts isolate each eigenvalue, and the recurrence
refines its eigenvector. The certified residual tends to zero while the
non-target distance remains positive. Matrix and integration errors
also tend to zero. Thus $\rho/g\to0$ and $B\to0$, whereas
$\widehat d\to d>0$. Both stopping tests eventually hold.
For unit vectors with acute angle $\theta$, the norm distance after
sign choice is $2\sin(\theta/2)\le\sqrt2\sin\theta$.
The alternating map preserves this distance by \eqref{eq:norm}.
\end{proof}

This is an input-dependent finite-dimensional assertion. It gives no
uniform working precision as $N$ grows. Tied concentration values do
not affect the scalar eigenpair certificates. Taking maxima or minima
of the respective enclosure endpoints also encloses the extremal
value without identifying a unique extremal index.

For \eqref{eq:weights}, write
$d\nu_R(t)=d\nu(t)$, $d\nu_S(t)=t\,d\nu(t)$ and let $p_i$ be the normalized
scalar eigenpolynomial with coefficient vector $v_i$ in the
$\widehat R$ or $\widehat S$ basis. Then
$$
\lambda_i=\int_0^s|p_i(t)|^2\,d\nu_j(t),\quad
1-\lambda_i=\int_s^b|p_i(t)|^2\,d\nu_j(t),\quad j=R,S.
\eqtag{eq:positive-eigenvalues}
$$
Each integral has a positive $N_j$-node quadrature rule exact through
degree $2N_j-1$ \cite[Section~3.5(v)]{DLMF}, whereas
$\deg p_i^2\le2N_j-2$.
Summing \eqref{eq:positive-eigenvalues} over an orthonormal eigenbasis gives
$\operatorname{tr}\op{G}_j$ and $\operatorname{tr}(\Id-\op{G}_j)$.
To compute these integrals without restricted-weight moment matrices,
put $\beta_R=\beta$, $\beta_S=\beta+1$. For Jacobi input and a scalar
polynomial $p$, use
$$
\int_0^s p(t)^2t^{\beta_j}(1-t)^\alpha\,dt
=s^{\beta_j+1}\int_0^1u^{\beta_j}(1-su)^\alpha p(su)^2\,du.
\eqtag{eq:jacobi-positive-rule}
$$
For the weight $u^{\beta_j}$, the nodes are the eigenvalues of its
tridiagonal recurrence matrix. If $\alpha=d$ is a nonnegative integer, $m$ nodes are exact
when $2m-1\ge2\deg p+d$. For Laguerre deficits, use
$$
\int_s^\infty p(t)^2t^{\beta_j}e^{-t}\,dt
=e^{-s}\int_0^\infty e^{-u}(s+u)^{\beta_j}p(s+u)^2\,du
\eqtag{eq:laguerre-positive-rule}
$$
and quadrature for the weight $e^{-u}$. The weights are $w_\ell^{-1}=\sum_{k=0}^{m-1}\pi_k(u_\ell)^2$, where
$(\pi_k)$ is orthonormal for the base quadrature weight
\cite[Sections~3.5(v)--(vi)]{DLMF}.
The error $e_Q$ includes the quadrature remainder when the rule is not
exact, as well as node, weight and recurrence-evaluation errors.
For Meixner, the interior integral is the finite positive sum in
Example~\ref{ex:meixner}. Deficits are integrated directly, without
subtracting a rounded value from one.
\begin{proposition}\label{prop:noninteger-quadrature}
Let $a,c>-1$, $0<h<1$, and let $p$ be a nonzero real polynomial of
degree at most $d$. Set
$$
J(p)=\int_0^1u^c(1-hu)^ap(u)^2\,du,\quad
J_0(p)=\int_0^1u^cp(u)^2\,du.
$$
Let $(u_\ell,w_\ell)_{\ell=1}^m$, with $0<u_\ell<1$, be the
positive quadrature rule for $u^c\,du$ exact through degree $2m-1$, and put
$$
Q_m(p)=\sum_{\ell=1}^m w_\ell(1-hu_\ell)^ap(u_\ell)^2.
$$
Choose an integer $K\ge0$ with $K+1>a$ and $2m-1\ge K+2d$, and define
$$
E_K=\frac{|(-a)_{K+1}|h^{K+1}}{(K+1)!(1-h)},\quad
\ell_a=\min\{1,(1-h)^a\}.
$$
Then $J_0(p)=\sum_\ell w_\ell p(u_\ell)^2$ and
$$
|J(p)-Q_m(p)|\le2E_KJ_0(p),\quad
\left|\frac{Q_m(p)}{J(p)}-1\right|\le\frac{2E_K}{\ell_a}.
\eqtag{eq:noninteger-quadrature-bound}
$$
If $a$ is a nonnegative integer and $K\ge a$, the rule is exact.
\end{proposition}

\begin{proof}
In the binomial series $(1-hu)^a=\sum_{k\ge0}(-a)_k(hu)^k/k!$,
the absolute ratios of successive tail terms are at most $h$ once
$k\ge K+1>a$. The remainder after degree $K$ is therefore bounded
uniformly on $0\le u\le1$ by $E_K$.
The rule integrates the truncated factor times $p^2$ exactly.
Bounding the remainder in both the integral and the quadrature gives
$2E_KJ_0(p)$. Since $(1-hu)^a\ge\ell_a$, the relative bound follows.
\end{proof}

For Jacobi input, the proposition applies to
\eqref{eq:jacobi-positive-rule}. Deficits have the direct form
$$
\begin{aligned}
\int_s^1p(t)^2t^{\beta_j}(1-t)^\alpha\,dt
&=(1-s)^{\alpha+1}\int_0^1u^\alpha
\bigl(1-(1-s)u\bigr)^{\beta_j}
 p\bigl(1-(1-s)u\bigr)^2\,du.
\end{aligned}
\eqtag{eq:jacobi-deficit-rule}
$$
The same bound applies with $h=1-s$, $a=\beta_j$ and $c=\alpha$.
It is independent of the polynomial coefficients and the magnitude
of the positive form. Node, weight and evaluation errors must still
be added to obtain $e_Q$ in \eqref{eq:accuracy-bound}.

\begin{example}\label{ex:collision-stability}
Take $N=3$, $\mu=4$ and $s=2+\varepsilon$, $0<\varepsilon<3/2$.
The unmodified Hermite matrix and its relevant gaps are
$$
\begin{gathered}
\op{L}=\begin{pmatrix}
9&0&3\sqrt2\\[7pt]0&3+\varepsilon&0\\[7pt]
3\sqrt2&0&6+2\varepsilon
\end{pmatrix},\quad
D=\sqrt{81-12\varepsilon+4\varepsilon^2},\\[7pt]
\delta=\frac{2\varepsilon(3-\varepsilon)}{9+D}
\sim\frac{\varepsilon}{3}.
\end{gathered}
$$
Here $D$ is the even-block gap and $\delta$ separates its lower
eigenvalue from the odd eigenvalue.
Let $v$ be the corresponding unit even vector and let $e$ be the odd
coordinate vector. For $\eta>0$, the symmetric perturbation
$\eta(ve^{\mathsf T}+ev^{\mathsf T})$ has norm $\eta$ and rotates the
nearby pair through an angle $\theta$ satisfying
$\tan(2\theta)=2\eta/\delta$, $0<\theta<\pi/4$.
At $\varepsilon=10^{-12}$ and $\eta=10^{-14}$, $\sin\theta>0.02995$.
Every parity-preserving perturbation of norm at most $\eta$, however,
has even-vector sine error at most $\eta/(D-\eta)<1.112\times10^{-15}$;
the odd vector is unchanged.

Choose $\zeta_R=(9+3\sqrt2)/2$, $\zeta_S=1$ in
\eqref{eq:normalized-commutant}. The normalized full gap is
$D/(9+3\sqrt2)>0.6796$ at these parameters. Even an arbitrary symmetric
perturbation of norm at most $\eta$ then gives sine error below
$1.472\times10^{-14}$.
These compare equal absolute perturbation norms. A full-matrix solver
that preserves parity also has the block protection; the example
quantifies sensitivity rather than a universal difference between solvers.
\end{example}

\begin{example}\label{ex:numerical-comparison}
Take the $N=12$ examples of \cite[Section~4]{CG}: Laguerre input
$(\beta,s)=(7/2,49)$ and Jacobi input
$(\alpha,\beta,s)=(1,1/2,1/256)$. We compare full diagonalization of
$\op{L}=\op{Q}^{\mathsf T}(\op{L}_R\oplus \op{L}_S)\op{Q}$ with its two
blocks \eqref{eq:exact}, and write $\op{M}=\op{M}_0(\chi)$.
The second block uses $\beta+1$.

After rounding 100-digit $\op{L}$ to binary64 and extracting its blocks,
both solves used LAPACK \texttt{dsyev} (Apple Accelerate, lower triangle).
Reference eigenvectors of the unrounded $\op{L}$, normalization, polynomial
bases and integration use 100 digits. Vectors are compared in
$\mathbb R^N$, with block vectors embedded using $\op{Q}^{\mathsf T}$.
Adaptive quadrature in \eqref{eq:positive-eigenvalues} was checked against
incomplete-gamma or incomplete-beta moments; all values were checked at
130 digits.

Let $v_i,\widehat v_i$ be the unit reference and computed eigenvectors,
matched by their eigenvalues of $\op{L}$. Set
$E_v=\max_i\sin\angle(\widehat v_i,v_i)$ and
$E_{\rm rel}=|\widehat d/d-1|$.
Here $d=1-\lambda_{\max}(\op{M})$ for Hermite and
$d=\lambda_{\min}(\op{M})$ for Gegenbauer, with reference values
$$
d_{\mathrm H}=1.09428890155\ldots\times10^{-29},\quad
d_{\mathrm G}=2.93714201786\ldots\times10^{-35}.
$$
We evaluate $\widehat d$ by the appropriate positive integral in
\eqref{eq:positive-eigenvalues}, at the computed vector matched to the
reference extremal eigenpair. For full-$\op{L}$ vectors, both scalar
components are included.
The errors are
$$
\begin{array}{llcc}
&\text{eigensolve}&E_v&E_{\rm rel}\\[7pt]
\text{Hermite}&\text{full }\op{L}
&1.76\times10^{-15}&2.03\times10^{-19}\\[7pt]
&\text{two parity blocks}
&4.05\times10^{-16}&1.23\times10^{-17}\\[7pt]
\text{Gegenbauer}&\text{full }\op{L}
&1.26\times10^{-15}&1.16\\[7pt]
&\text{two parity blocks}
&6.88\times10^{-16}&1.01\times10^{-2}.
\end{array}
$$
The separate-block solves give smaller $E_v$ in both tests. For the
smallest Gegenbauer eigenvalue, they reduce $E_{\rm rel}$ from $116\%$
to $1.01\%$; the full solve gives the smaller relative error in the
Hermite deficit.
\end{example}

\begin{samepage}
\begin{example}\label{ex:validated-computation}
For Jacobi input $(\alpha,\beta)=(1/2,1/4)$ and scalar order $14$,
direct binary64 Cholesky factorization of the monomial moment matrix
fails. The basis defined by binary64 recurrence data has Frobenius
orthogonality error $3.52\times10^{-15}$, checked with 100-digit inner
products.

For the $N=12$ Gegenbauer example above, Algorithm~\ref{alg:manufacture}
computed the scalar eigenvectors at 30 decimal digits and verified the
residuals and positive forms by outward-rounded interval arithmetic at
90 digits. The nodes were enclosed by sign counts;
weights were evaluated by the reciprocal sums above.
Since $\alpha=1$ and $\deg p\le5$, the integrand in
\eqref{eq:jacobi-positive-rule}, apart from $u^{\beta_j}$, has degree
at most $11$. Six nodes per block therefore suffice exactly.
All twelve concentration values were enclosed, isolating the minimum.
For its computed value
$$
\widehat d=2.93714201785750961572874550537485\ldots\times10^{-35},
$$
the verified bounds, with $\Omega_S=1$, were
$$
\begin{aligned}
\rho&<1.080\times10^{-29},\quad g>11.476,\quad e_Q<2.20\times10^{-86},\\[7pt]
B&<8.86\times10^{-61},\quad
\left|\frac{\widehat d}{\lambda_{\min}(\op{M})}-1\right|
<3.02\times10^{-26}.
\end{aligned}
$$
This computation uses higher precision than the preceding binary64
eigensolves. The tolerance is verified by \eqref{eq:accuracy-bound};
the alternating reconstruction transports the enclosure to every $\tau$.
It does not change the intrinsic conditioning of the scalar problem.
\end{example}
\end{samepage}

\begin{example}\label{ex:noninteger-stress}
For $N=24$, we tested noninteger Jacobi parameters near $-1$ and
cutoffs near the endpoints. Scalar eigenvectors were computed at
50 decimal digits; outward interval arithmetic at 110 digits
certified their residuals, all within-block eigenvalue separations
and the positive forms. Proposition~\ref{prop:noninteger-quadrature} bounds the
quadrature remainder. Every concentration value, or deficit, was
enclosed before isolating the minimum.

The first two rows give $d=\lambda_{\min}(\op{M})$ and the last gives
$d=1-\lambda_{\max}(\op{M})$. Leading digits of $\widehat d$ are shown;
the relative bounds apply to the complete computed values.
$$
\begin{array}{ccccc}
(\alpha,\beta)&s&\widehat d&(m_R,m_S)&|\widehat d/d-1|\\[7pt]
(-0.99,-0.999)&1/4&3.118051652\ldots\times10^{-28}
 &(66,66)&<1.804\times10^{-53}\\[7pt]
(1/3,-3/4)&1/256&4.959531999\ldots\times10^{-68}
 &(25,25)&<8.847\times10^{-32}\\[7pt]
(-0.99,-0.999)&255/256&3.725544103\ldots\times10^{-69}
 &(25,24)&<4.327\times10^{-29}.
\end{array}
$$
Here $m_R,m_S$ are quadrature-node counts. The truncation degrees
$K$ were $(108,108)$, $(26,26)$ and $(27,25)$, respectively.
An independent calculation at 165 digits using incomplete-beta
moments lay inside each certified enclosure.
The reconstruction transports these bounds to every $\tau$.
\end{example}

\begin{example}\label{ex:large-dimensions}
The normalized construction was also tested at larger dimensions for
the original Hermite parameters $\mu=4$, $s=49$, and Jacobi parameters
$\alpha=1$, $\beta=1/2$, $s=1/256$.
Here both routes diagonalize the same
$\widehat{\op{L}}$ from \eqref{eq:normalized-commutant}:
the full route uses the symmetric band driver \texttt{DSBEV} with
half-bandwidth two, and the scalar route uses two calls to the
tridiagonal driver \texttt{DSTEV}, all in binary64.
Outward interval arithmetic at 80 decimal digits, applied to the exact
recurrence coefficients, encloses the matrix construction error,
the residual of every returned vector, and all scalar eigenvalues
using sign counts. Let $v_i$ be the exact unit eigenvectors of
$\widehat{\op{L}}$ and $\widehat v_i$ the
normalized computed vectors with the corresponding ordered labels.
The table gives upper bounds for the matrix norm and eigenvector sine
error, and lower bounds for the smallest gap.
$$
\begin{array}{ccccc}
\text{family}&N&\|\widehat{\op{L}}\|_2
 &\operatorname{gap}(\widehat{\op{L}})
 &\max_i\sin\angle(\widehat v_i,v_i)\\[7pt]
\text{Hermite}&64 &<2.985174&>1.4767\times10^{-2}&<8.824\times10^{-14}\\[7pt]
&128&<2.989508&>7.0061\times10^{-3}&<2.172\times10^{-13}\\[7pt]
&256&<2.993666&>2.9084\times10^{-3}&<7.447\times10^{-13}\\[7pt]
\text{Jacobi}&64 &<2.967768&>2.7336\times10^{-3}&<3.863\times10^{-13}\\[7pt]
&128&<2.983079&>1.1707\times10^{-3}&<1.013\times10^{-12}\\[7pt]
&256&<2.991344&>5.1841\times10^{-4}&<3.197\times10^{-12}.
\end{array}
$$
The two drivers returned identical binary64 eigenpairs, up to signs,
on these six inputs, so these tests show no accuracy advantage for the
scalar route over an efficient banded driver.
The scalar representation uses $2N-2$ diagonal and off-diagonal
entries, compared with $3N$ allocated entries in the band format;
retaining the two eigenvector matrices in scalar form uses
$N_R^2+N_S^2$ entries, compared with $N^2$ for the full eigenvector
matrix.  The certificates concern eigenvectors and gaps; no relative
accuracy for the extreme concentration values is asserted at these
larger dimensions. By \eqref{eq:norm}, the angle bounds also hold for
the reconstructed unit eigenpolynomials for every $\tau$.

\begin{samepage}
At $N=256$, median times over nine runs on Apple Accelerate were
$4.36$ versus $2.21$ milliseconds for Hermite and $4.06$ versus $2.18$
for Jacobi, for the full and block routes respectively. Both wrappers
include allocation and copies; the block route also interlaces its
vectors. At $N=64$ the timings were essentially equal. These timings
exclude recurrence construction and interval verification, and are
specific to the inputs, drivers and platform.

\end{samepage}
\end{example}

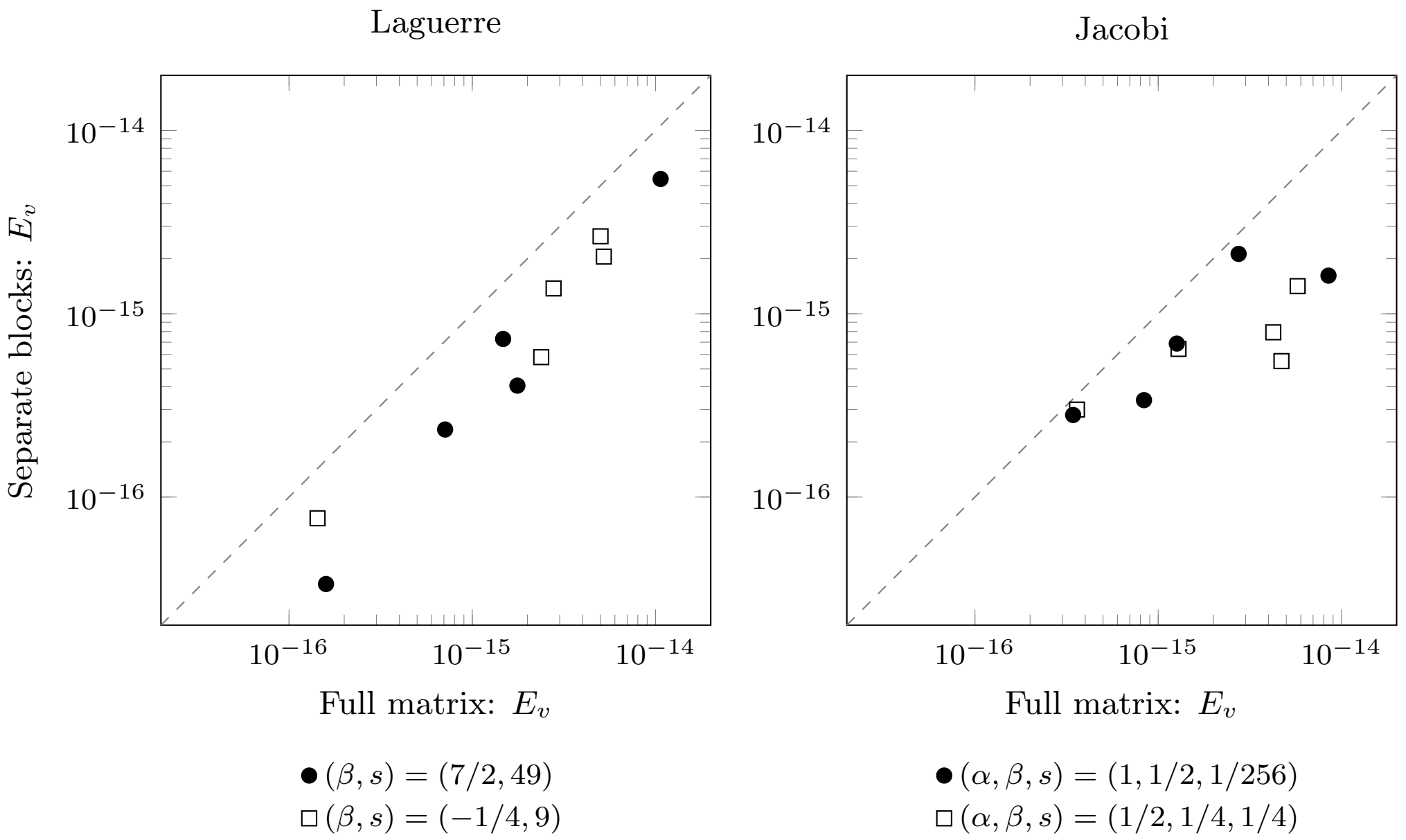
\begin{figure}[htbp]
\centering
\begin{tikzpicture}[font=\footnotesize]
\begin{groupplot}[
group style={group size=2 by 1,horizontal sep=1.3cm},
scale only axis,width=0.35\textwidth,height=0.35\textwidth,
xmode=log,ymode=log,
xmin=2e-17,xmax=2e-14,ymin=2e-17,ymax=2e-14,
xtick={1e-16,1e-15,1e-14},ytick={1e-16,1e-15,1e-14},
minor tick num=0,
xlabel={Full matrix: $E_v$},
tick label style={font=\scriptsize},
legend style={at={(0.5,-0.22)},anchor=north,draw=none,
font=\scriptsize,legend cell align=left},
]
\nextgroupplot[title={Laguerre},ylabel={Separate blocks: $E_v$}]
\addplot[gray,dashed,forget plot] coordinates {(2e-17,2e-17) (2e-14,2e-14)};
\addplot[black,only marks,mark=*,mark size=2pt] coordinates
{(1.59181051e-16,3.35232273e-17) (7.10772890e-16,2.33644229e-16)
(1.76299997e-15,4.05441767e-16) (1.47287259e-15,7.28960697e-16)
(1.06583199e-14,5.44085322e-15)};
\addlegendentry{$(\beta,s)=(7/2,49)$}
\addplot[black,only marks,mark=square,mark size=2pt] coordinates
{(1.43130833e-16,7.66060449e-17) (2.37725202e-15,5.79885259e-16)
(2.77687676e-15,1.37516984e-15) (5.22145836e-15,2.05099530e-15)
(5.00971804e-15,2.64740438e-15)};
\addlegendentry{$(\beta,s)=(-1/4,9)$}
\nextgroupplot[title={Jacobi}]
\addplot[gray,dashed,forget plot] coordinates {(2e-17,2e-17) (2e-14,2e-14)};
\addplot[black,only marks,mark=*,mark size=2pt] coordinates
{(3.43502831e-16,2.80382970e-16) (8.36574398e-16,3.37889271e-16)
(1.26347721e-15,6.87956996e-16) (2.74710807e-15,2.12199849e-15)
(8.49738638e-15,1.61567740e-15)};
\addlegendentry{$(\alpha,\beta,s)=(1,1/2,1/256)$}
\addplot[black,only marks,mark=square,mark size=2pt] coordinates
{(3.60553506e-16,3.00152499e-16) (1.29072499e-15,6.42941858e-16)
(4.70516610e-15,5.51920136e-16) (4.24588937e-15,7.92360309e-16)
(5.76673348e-15,1.41570430e-15)};
\addlegendentry{$(\alpha,\beta,s)=(1/2,1/4,1/4)$}
\end{groupplot}
\end{tikzpicture}
\caption{Eigenvector errors computed here for full (horizontal) and
separate-block (vertical) diagonalizations of the commuting matrices
in \cite{CG}. Both use identical binary64 data and the same QR driver.
Points below the diagonal have smaller $E_v$ for the block solves.}
\label{fig:comparison}
\end{figure}

Using the binary64 protocol of Example~\ref{ex:numerical-comparison},
Figure~\ref{fig:comparison} compares four inputs at
$N=6,9,12,17,24$. The blocks give smaller $E_v$ in all twenty tests,
by factors between $1.20$ and $8.53$; their largest error is
$5.44\times10^{-15}$.

The ordering of $E_v$ depends on the input. With the same protocol,
$N=5$ and Jacobi input $(\alpha,\beta,s)=(1/2,1/4,1/4)$, the full solve
gives $E_v=2.85\times10^{-16}$, compared with $3.00\times10^{-16}$
for the block solves. These values were checked at 100 and 130 digits.

\begin{samepage}
For the Jacobi input $(\alpha,\beta,s)=(1/2,1/4,1/4)$,
$N=9$ and $\tau=1$, we also integrated the extremals of the two
scalar blocks, reconstructed in $x$. The density is positive on
$(-\sqrt2,-1)\cup(1,\sqrt2)$, and the window is $(-\sqrt5/2,-1)\cup(1,\sqrt5/2)$.
The largest block eigenvalues are
$$
\begin{aligned}
\lambda_R&=\lambda_{\max}(\op{G}_R)=0.997377870710549\ldots,\\[7pt]
\lambda_S&=\lambda_{\max}(\op{G}_S)=0.917749195280867\ldots.
\end{aligned}
$$
At 100 digits, the reconstructed norms and retained energies agree
with $1$ and these values within $10^{-95}$.

\end{samepage}

\noindent\begin{minipage}{\textwidth}
\section*{Acknowledgements}
The author is deeply grateful to one of his thesis supervisors, Kenier
Castillo, for teaching him to seek simplicity and generality through duality.
Castillo's handwritten notes on this problem already contained the central
ideas developed here. He declined coauthorship.

The author was supported by FAPESP, grants 2024/01546-1 and 2024/23674-1.

\end{minipage}

\begin{samepage}
\section*{Declaration on the Use of Artificial Intelligence}
ChatGPT assisted with writing. Responsibility for the final manuscript rests with the
author. Further extensions may be explored with AI when a specific problem calls for
them. Our aim is to use the alternating map to identify the common structure behind a
broad range of problems.

\end{samepage}


\begin{thebibliography}{99}
\bibitem{Admissible}
K.~Castillo,
\emph{Half-step-invariant sets, admissible maps, and classical orthogonal polynomials},
preprint, arXiv:2605.27542v1, May~2026, 77~pp.

\bibitem{Manufacture}
K.~Castillo and G.~Gordillo-N\'u\~nez,
\emph{The manufacture of examples: the $(-1)$-classical orthogonal polynomials},
Res. Math. Sci. \textbf{13} (2026), Paper No.~81.

\bibitem{CG}
M.~M.~Castro and F.~A.~Gr\"unbaum,
\emph{Time-and-band limiting for symmetric classical Dunkl orthogonal polynomials},
preprint, arXiv:2610.24100v1, 2026.

\bibitem{DavisKahan}
C.~Davis and W.~M.~Kahan,
\emph{The rotation of eigenvectors by a perturbation. III},
SIAM J. Numer. Anal. \textbf{7} (1970), no.~1, 1--46.

\bibitem{PolynomialMappings}
M.~N.~de Jesus and J.~Petronilho,
\emph{On orthogonal polynomials obtained via polynomial mappings},
J. Approx. Theory \textbf{162} (2010), no.~12, 2243--2277.

\bibitem{GVZChihara}
V.~X.~Genest, L.~Vinet, and A.~Zhedanov,
\emph{A \textquotedblleft{}continuous\textquotedblright{} limit of the complementary Bannai--Ito polynomials:
Chihara polynomials},
SIGMA \textbf{10} (2014), Paper No.~038, 18~pp.

\bibitem{Grunbaum}
F.~A.~Gr\"unbaum,
\emph{A new property of reproducing kernels for classical orthogonal polynomials},
J. Math. Anal. Appl. \textbf{95} (1983), no.~2, 491--500.

\bibitem{GVZ}
F.~A.~Gr\"unbaum, L.~Vinet, and A.~Zhedanov,
\emph{Algebraic Heun operator and band-time limiting},
Comm. Math. Phys. \textbf{364} (2018), no.~3, 1041--1068.

\bibitem{Jafarov}
E.~I.~Jafarov and J.~Van der Jeugt,
\emph{Discrete series representations for $\mathfrak{sl}(2|1)$,
Meixner polynomials and oscillator models},
J. Phys. A: Math. Theor. \textbf{45} (2012), Paper No.~485201.

\bibitem{Maroni1993}
P.~Maroni,
\emph{Variations around classical orthogonal polynomials. Connected problems},
J. Comput. Appl. Math. \textbf{48} (1993), no.~1--2, 133--155.

\bibitem{DLMF}
\emph{NIST Digital Library of Mathematical Functions},
National Institute of Standards and Technology,
Sections~3.5(v)--(vi) and~18.22(ii).

\bibitem{Perlstadt}
M.~Perlstadt,
\emph{Chopped orthogonal polynomial expansions---some discrete cases},
SIAM J. Algebraic Discrete Methods \textbf{4} (1983), no.~1, 94--100.
\end{thebibliography}
\end{document}